\documentclass{article}

\usepackage{amssymb,amscd,amsthm,verbatim,fancyhdr}
\usepackage[leqno]{amsmath}
\usepackage{mathrsfs}

\begin{document}
\numberwithin{equation}{section}
\newtheorem{thm}{Theorem}[section]
\newtheorem{lemma}[thm]{Lemma}
\newtheorem{clm}[thm]{Claim}
\newtheorem{remark}[thm]{Remark}
\newtheorem{definition}[thm]{Definition}
\newtheorem{cor}[thm]{Corollary}


\newcommand{\hdt}{{\dot{\mathrm{H}}^{1/2}}}
\newcommand{\hdtr}{{\dot{\mathrm{H}}^{1/2}(\mathbb{R}^3)}}
\newcommand{\R}{\mathbb{R}}
\newcommand{\ei}{\mathrm{e}^{it\Delta}}
\newcommand{\ltrt}{L^3(\mathbb{R}^3)}
\newcommand{\lt}{L^3}
\newcommand{\rt}{\mathbb{R}^3}
\newcommand{\X}{\mathfrak{X}}
\newcommand{\F}{\mathfrak{F}}
\newcommand{\hdhalf}{{\dot H^\frac{1}{2}}}
\newcommand{\hdthalf}{{\dot H^\frac{3}{2}}}
\newcommand{\hdo}{\dot H^1}
\newcommand{\rthmiz}{\R^3\times(-\infty,0)}
\newcommand{\q}[2]{{#1}_{#2}}
\renewcommand{\t}{\theta}
\newcommand{\lxt}[2]{L_{x,\,t}^{#1}}
\newcommand{\rr}{\sqrt{x_1^2+x_2^2}}
\newcommand{\ve}{\varepsilon}
\newcommand{\hdhrt}{\dot H^\frac{1}{2}(\mathbb{R}^3)}
\renewcommand{\P}{\mathbb{P} }
\newcommand{\RR}{\mathcal{R} }
\newcommand{\TT}{\overline{T} }
\newcommand{\e}{\epsilon }
\newcommand{\D}{\Delta }
\renewcommand{\d}{\delta }
\renewcommand{\l}{\lambda }
\newcommand{\To}{\TT_1 }
\newcommand{\ukt}{u^{(KT)} }
\newcommand{\ttil}{\tilde{T_1} }
\newcommand{\etl}{e^{t\Delta} }
\newcommand{\et}{\mathscr{E}_T }
\newcommand{\se}{\mathscr{E}}
\newcommand{\ft}{\mathscr{F}_T }
\newcommand{\eti}{\mathscr{E}^{\infty}_T }
\newcommand{\fti}{\mathscr{F}^{\infty}_T }
\newcommand{\xjn}{x_{j,n}}
\newcommand{\xjpn}{x_{j',n}}
\newcommand{\ljn}{\l_{j,n} }
\newcommand{\ljpn}{\l_{j',n} }
\newcommand{\lkn}{\l_{k,n} }
\newcommand{\voj}{V_{0,j} }
\newcommand{\voo}{V_{0,1} }
\newcommand{\uon}{u_{0,n} }
\newcommand{\N}{\mathbb{N} }
\newcommand{\E}{\mathscr{E} }
\newcommand{\tu}{\tilde{u} }
\newcommand{\tU}{\tilde{U} }
\newcommand{\etj}{{\E_{T^*_j}}}
\newcommand{\B}{\mathcal{B}}
\newcommand{\tujn}{\tilde{U}_{j,n}}
\newcommand{\tujpn}{\tilde{U}_{j',n}}
\newcommand{\soj}{\sum_{j=1}^J}
\newcommand{\soi}{\sum_{j=1}^\infty}
\newcommand{\hnj}{H_{n,J}}
\newcommand{\enj}{e_{n,J}}
\newcommand{\pnj}{p_{n,J}}
\newcommand{\lfoi}{L^5_{(0,\infty)}}
\newcommand{\lfhoi}{L^{{5/2}}_{(0,\infty)}}
\newcommand{\wnj}{w_n^J}
\newcommand{\rnj}{r_n^J}
\newcommand{\lfij}{{L^5_{I_j}}}
\newcommand{\lfiz}{{L^5_{I_0}}}
\newcommand{\lfio}{{L^5_{I_1}}}
\newcommand{\lfit}{{L^5_{I_2}}}
\newcommand{\lfhij}{{L^{{5/2}}_{I_j}}}
\newcommand{\lfhiz}{{L^{{5/2}}_{I_0}}}
\newcommand{\lfhio}{{L^{{5/2}}_{I_1}}}
\newcommand{\lfhit}{{L^{{5/2}}_{I_2}}}
\newcommand{\lfi}{{L^5_{I}}}
\newcommand{\lfhi}{{L^{{5/2}}_{I}}}
\newcommand{\doh}{D^\frac{1}{2}}
\newcommand{\tjn}{t_{j,n}}
\newcommand{\tjpn}{t_{j',n}}
\newcommand{\wnlj}{w_n^{l,J}}
\renewcommand{\O}{\mathcal{O}}
\newcommand{\tn}{{\tilde \|}}
\newcommand{\xto}{\xrightarrow[n\to\infty]{}}

\title{An alternative approach to regularity for the Navier-Stokes equations in critical spaces}
\author{Carlos E. Kenig\footnote{Department of Mathematics, University of Chicago; Chicago, Il 60637, USA; cek@math.uchicago.edu; supported in part by NSF grant DMS-0456583} \and
Gabriel S. Koch\footnote{Department of Mathematics, University of Chicago; Chicago, Il 60637, USA; koch@math.uchicago.edu}}

\date{}
\maketitle
\begin{abstract}
In this paper we present an alternative viewpoint on recent studies of regularity of solutions to the Navier-Stokes equations in critical spaces.  In particular, we prove that mild solutions which remain bounded in the space $\hdhalf$ do not become singular in finite time, a result which was proved in a more general setting by L. Escauriaza, G. Seregin and V. \v Sver\' ak using a different approach. We use the method of ``concentration-compactness" $+$ ``rigidity theorem" which was recently developed by C. Kenig and F. Merle to treat critical dispersive equations.  To the authors' knowledge, this is the first instance in which this method has been applied to a parabolic equation.
\\\\\\
\textbf{R\'esum\'e:} \quad Dans cet expos\'e, nous pr\'esentons un point de vue diff\'erent sur les \'etudes r\'ecentes concernant la r\'egularit\'e des solutions des \'equations de Navier-Stokes dans les espaces critiques.  En particulier, nous d\'emontrons que les solutions faibles qui restent born\'ees dans l'espace $\dot H^{1/2}$ ne deviennent pas singuli\`eres en temps fini. Ce r\'esultat a \'et\'e d\'emontr\'e dans un cas plus g\'en\'eral par L. Escauriaza, G. Seregin et V. \v Sver\' ak en utilisant une approche diff\'erente.  Nous utilisons la m\'ethode de $\ll$ concentration-compacit\'e $\gg$ $+$ $\ll$ th\'eor\`eme de rigidit\'e $\gg$ qui a \'et\'e r\'ecemment d\'evelopp\'ee par C. Kenig et F. Merle pour traiter les \'equations dispersives critiques.  \`A la connaissance des auteurs, c'est la premi\`ere fois que cette m\'ethode est appliqu\'ee \`a une \'equation parabolique.
\end{abstract}
\section*{Introduction}
In recent studies, the idea of establishing the existence of so-called ``critical elements" (or the earlier ``minimal blow-up solutions") has led to significant progress in the theory of ``critical" dispersive and hyperbolic equations such as the energy-critical nonlinear Schr\" odinger equation \cite{bourgain1,bourgain2,ckstt1,km1,rv,visan}, mass-critical nonlinear Schr\" odinger equation \cite{keraani2,ktv,kvz,tvz2,tvz}, $\hdhalf$-critical nonlinear Schr\" odinger equation \cite{km3}, energy-critical nonlinear wave equation \cite{km2}, energy-critical and mass-critical Hartree equations \cite{mxz1,mxz2,mxz3,mxz4,mxz5} and energy-critical wave maps \cite{ckm,ks,tao1,tao2}.
\\\\
In this paper we exhibit the generality of the method of ``critical elements" by applying it to a parabolic system, namely the standard\footnote{For simplicity, we have set the coefficient of kinematic viscosity $\nu =1$.} Navier-Stokes equations (NSE):
\begin{equation}\label{eq:nseA}
\begin{array}{rcl}
u_t -\D u + (u\cdot \nabla ) u+\nabla p & = &  0 \\
\nabla \cdot u & = & 0
\end{array}
\end{equation}
Typically, $u$ is interpreted as the velocity vector field of a fluid filling a region in space, and $p$ is the associated scalar pressure function.
\\\\
The ``critical" spaces are those which are invariant under the natural scaling of the equation.  For NSE, if $u(x,t)$ is a solution, then so is $u_\lambda(x,t):= \lambda u(\lambda x, \lambda^2 t)$ for any $\lambda >0$.  The critical spaces are of the type $\X$, where $\|u_\lambda\|_\X = \|u\|_\X$.  For the Navier-Stokes equations, one can take, for example, $\X=L^\infty((0,+\infty);L^3(\mathbb{R}^3))$ or the smaller space $\X=L^\infty((0,+\infty); \dot H^\frac{1}{2}(\mathbb{R}^3))$.  In fact, one has the ``chain of critical spaces" given by the continuous embeddings
$$\hdhalf(\rt) \hookrightarrow \ltrt \hookrightarrow  \dot B^{-1+\frac{3}{p}}_{p,\infty}(\rt)_{(p<\infty)} \hookrightarrow BMO^{-1}(\rt) \hookrightarrow \dot B^{-1}_{\infty,\infty}(\rt)\ .$$
These are the spaces in which the initial data of solutions in the critical settings live, and we will also refer to them as ``critical spaces" -- that is, spaces of functions on $\rt$ whose norms satisfy $\|\l f(\l \cdot)\| = \|f\|$.
\\\\
In the recent important paper \cite{ess}, L. Escauriaza, G. Seregin and V. \v Sver\' ak showed that any ``Leray-Hopf" weak solution which remains bounded in $L^3(\mathbb{R}^3)$ cannot develop a singularity in finite time.  Their proof used a blow-up procedure and reduction to a backwards uniqueness question for the heat equation, and was then completed using Carleman-type inequalities and the theory of unique continuation.  Here, we approach the same problem using the method of ``critical elements".  Although we do also use the main tools appearing in \cite{ess} to complete our proof, we hope that it is more intuitively clear in our exposition why those particular tools are needed.
\\\\
The precise statement of the main result we address in this paper is the following:
\begin{thm}\label{thm:regularity}
Let $u_0 \in \hdhalf$ satisfy $\nabla \cdot u_0 = 0$ and let $u=NS(u_0)$ be the associated ``mild solution'' to NSE satisfying $u(0)=u_0$.  Suppose that there is some $A>0$ such that $\|u(t)\|_{\hdhalf(\rt)} \leq A$ for all $t>0$ such that $u$ is defined.  Then $u$ is defined (and smooth) for all positive times.
\end{thm}
\noindent
(Theorem \ref{thm:regularity} of course follows from the result in \cite{ess}.)  We believe the methods given below will work as well if we replace $\hdhalf$ by $\lt$ in Theorem \ref{thm:regularity} and hence can be used to give an alternative proof of the result in \cite{ess} in the case of mild solutions.  For technical reasons (described below) we start with the above result and plan to return to the more general case in a future publication.
\\\\
The ``mild solutions" to NSE which we consider in our approach have the form
\begin{equation}\label{eem}
u(t)= L(t)u_0 + \int_0^t L(t-s)f(u(s))\ ds
\end{equation}
for some divergence-free initial datum $u_0$, with the linear solution operator $L$ and nonlinearity $f$ given by
\begin{equation}\label{ns1}
L(t)= e^{t\triangle}\ , \qquad f(u(s)) = -\mathbb{P}(u\cdot \nabla u)(s)\ .
\end{equation}
Here $e^{t\triangle}u_0$ is the convolution of $u_0$ with the heat kernel, and the equation (\ref{eem}), (\ref{ns1}) comes formally from applying the Helmholtz projection operator $\P$ to (\ref{eq:nseA}) which fixes $u_t - \D u$ and eliminates the term $\nabla p$, and then solving the resulting non-homogeneous heat equation by Duhamel's formula.  Under sufficient regularity assumptions, mild solutions are in fact classical solutions, and the existence of such a solution on some time interval is typically established via the contraction-mapping principle in an appropriate function space.
\\\\
We follow the method of C. Kenig and F. Merle.  In a series of recent works \cite{km3,km1,km2}, they use the method of ``critical elements" to approach the question of global existence and ``scattering" (approaching a linear solution at large times) for nonlinear hyperbolic and dispersive equations in  critical settings.  For example, for the 3D non-linear Schr\" odinger equation $iu_t + \triangle u = f(u)$ (NLS), they considered (\ref{eem}) with
\begin{equation}\label{nls1}
L(t)=e^{it\triangle}\ ,
\end{equation}
the free Schr\" odinger operator, and both
\begin{equation}\label{nls2}
f(u) = -|u|^4u\  \qquad (\ \ \X=L^\infty((0,+\infty); \dot H^1(\mathbb{R}^3))\ \ )
\end{equation}
(the focusing case), and
\begin{equation}\label{nls3}
f(u) = |u|^2u\  \qquad (\ \ \X=L^\infty((0,+\infty); \dot H^\frac{1}{2}(\mathbb{R}^3))\ \ )
\end{equation}
(the defocusing case), the exponents needed for the ``$\dot H^1$-critical" and ``$\dot H^\frac{1}{2}$-critical" settings, respectively.  Note that in the case (\ref{nls3}) of cubic nonlinearity, the equation is invariant under the same scaling as the Navier-Stokes equations.
\\\\
In all cases, the general strategy was essentially the same, which we'll describe now in the $\hdhalf$ setting:
\\\\
For any $u_0\in \dot H^\frac{1}{2}(\mathbb{R}^3)$, one uses fixed-point arguments to assign a maximal time $T^*(u_0) \leq +\infty$ such that a solution $u$ to (\ref{eem}) which remains in $\hdhalf$ for positive time exists and is unique in some scaling-invariant space $\mathcal{X}_T$ for any fixed $T<T^*(u_0)$, where $\mathcal{X}_T$ denotes a space of functions defined on the space-time region $\rt \times (0,T)$.  Define
$$|||u|||:= \sup_{t\in[0,T^{*}(u_0))}\|u(t)\|_{\dot{H}^\frac{1}{2}}\ .$$
The type of result proved in \cite{km3} (variants of which were proved in \cite{km1,km2} and which we will prove here as well) is that $|||u||| < +\infty$ implies that $T^*(u_0) = +\infty$ and $\|u(t) - L(t)u_0^+\|_{\dot H^\frac{1}{2}(\mathbb{R}^3)} \to 0 $ as $t\to +\infty$ for some $u_0^+ \in \dot H^\frac{1}{2}(\mathbb{R}^3)$.  In other words, $u$ exists globally and scatters.  Typically this is known to be true for $|||u|||<\epsilon_0$ for sufficiently small $\epsilon_0 > 0 $.  In the case of NSE, the scattering condition is replaced by decay to zero in norm -- in other words, we set $u_0^+ = 0$.  For globally defined $\hdhalf$-valued solutions, such decay was proved in \cite{GIP1} (see also \cite{GIP2,GIP3}).
\\\\
It is worth pointing out that the $\hdhalf$ decay, which is proved quite easily in \cite{GIP1} by decomposing the initial data into a large part in $L^2$ and a small part in $\hdhalf$, solving the corresponding equations and employing the standard energy arguments, is actually used heavily in our proofs, and significantly reduces the difficulty from the NLS case treated in \cite{km3}.
\\\\
The general method of proof, which can be referred to as ``concentration-compactness'' $+$ ``rigidity'', is comprised of the following three main steps for a proof by contradiction:
\\
\begin{enumerate}
\item  \underline{Existence of a ``critical element"}\\\\ Assuming a finite maximal threshold $A_c >0$ for which $|||u||| < A_c$ implies global existence and scattering but such a statement fails for any $A>A_c$, there exists a solution $u_c$ with $|||u_c|||=A_c$, for which global existence or scattering fails.

\item  \underline{Compactness of critical elements} \\\\ Such a critical element $u_c$ produces a ``compact family" -- that is, up to norm-invariant rescalings and translations in space (and possibly in time), the set $\{u_c(t)\}$ is pre-compact in a critical space.

\item  \underline{Rigidity}\\\\ The existence of the compact family produces a contradiction to known results.\\
\end{enumerate}
Steps 1-2 are accomplished by considering a minimizing sequence of solutions $\{u_n\}$ with initial data $\{u_{0,n}\}$ such that $|||u_n|||\leq A_n$, $A_n \searrow A_c$ for which global existence or scattering fails (typically quantified by $\|u_n\|_{\mathcal{X}_{T^*}} = +\infty$).
\\\\
Then the main tool for realizing this program is a ``profile decomposition" associated to that sequence, which explores the lack of compactness in the embedding $\dot{H}^\frac{1}{2}(\mathbb{R}^3) \hookrightarrow L^3(\mathbb{R}^3)$.  For example, in \cite{km3} the following decomposition (based on \cite{keraani}) was used for treatment of the NLS case:  There exists a sequence
$\{V_{0,j}\}_{j=1}^\infty\subseteq\hdt$, with associated linear solutions $V_j^l(x,t)=\ei V_{0,j}$, and sequences of scales $\lambda_{j,n}\in \R^+$ and shifts $x_{j,n}\in \rt$ and $t_{j,n}\in\R$, such that (after a subsequence in $n$)
\begin{equation}\label{prof}
u_{0,n}(x)=\sum_{j=1}^J\frac{1}{\lambda_{j,n}}V_j^l\left(\frac{x-x_{j,n}}{\lambda_{j,n}},
-\frac{t_{j,n}}{\lambda_{j,n}^2}\right)+w_n^J(x)
\end{equation}
where the ``profiles" in the sum are orthogonal in a certain sense due to the choices of $\lambda_{j,n}$, $x_{j,n}$ and $t_{j,n}$, and $w_n^J$ is a small error for large $J$ and $n$.  (From this, a similar decomposition is then established for the time evolution of $\uon$.)  Using this setup, it is shown that in fact the solution to NLS with initial datum $V_{0,j_0}$ is a critical element for some $j_0 \in \N$, which establishes Step 1.  Compactness (Step 2) is established using the same tools.  (We will discuss Step 3 momentarily.)
\\\\
Our ultimate goal here is to give such a proof in the context of the Navier-Stokes equations, where we would consider (\ref{eem}), (\ref{ns1}) for $u_0 \in \ltrt$, with $\mathcal{X}_T=\mathcal{C}([0,T];\ltrt) \cap L^5(\rt \times (0,T))$ and $T^*(u_0)$ defined accordingly, and
$$|||u|||:= \sup_{t\in[0,T^*(u_0))}\|u(t)\|_{\ltrt}\ .$$
Since a profile decomposition has already been established by I. Gallagher in \cite{gallagher} for solutions to the Navier-Stokes equations evolving from a bounded set in $\hdhalf$ (based on the decomposition for the initial data in \cite{Gerard}), we restrict ourselves in this paper to the case of data in $\hdhalf$.  We expect that the result in \cite{gallagher} can be extended to the $L^3$ setting, which would then allow for an extension of our approach to that case.  We plan to return to this in a future publication\footnote{At the time of publication of this article, an $\ltrt$ profile decomposition has been established by the second author in \cite{gk}, and the program in $\ltrt$ has been completed in a collaboration of the second author with I. Gallagher and F. Planchon, to appear soon.}.
\\\\
We take  $$\mathcal{X}_T:=\mathcal{C}([0,T];\hdhalf(\rt)) \cap L^2((0,T);\hdthalf(\rt))$$ and as before let
$$|||u|||:= \sup_{t\in[0,T^*(u_0))}\|u(t)\|_{\hdhalf(\rt)}\ .$$
The profile decomposition of \cite{gallagher} for solutions corresponding to a bounded sequence $\{u_{0,n}\} \subset \hdhalf(\rt)$ of divergence-free fields takes the form (after a subsequence in $n$)
$$u_n(x,t) = \sum_{j=1}^J\frac{1}{\lambda_{j,n}}U_j\left(\frac{x-x_{j,n}}{\lambda_{j,n}},
\frac{t}{\lambda_{j,n}^2}\right)+e^{t\triangle}w_n^J(x) + r_n^J(x,t)\ ,$$
where $u_n$ and $U_j$ are the solutions to the Navier-Stokes equations with initial data $u_{0,n}$ and $V_{0,j}$, respectively, and $w_n^J$ and $r_n^J$ are again small errors for large $J$ and $n$.  (We remark that the absence in the above decomposition of the time shifts $\tjn$ which appeared in (\ref{prof}) greatly simplifies matters on a technical level; in \cite{km3}, this was a significant consideration.)  The above program is then completed in Theorems \ref{thm:a}, \ref{thm:b} and \ref{thm:c} as follows:
\\\\
By known local-existence and small-data results, there exists a small $\epsilon_0 >0$ such that $|||u||| < \epsilon_0$ implies that the solution exists for all time and tends to zero in the $\dot{H}^\frac{1}{2}$ norm.  We thereby assume a finite critical value $A_c \geq \epsilon_0 >0$ (as in Step 1) such that any solution $u$ with $|||u||| < A_c$ must exist globally and decay to zero in $\dot{H}^\frac{1}{2}$, and $A_c$ is the maximum such value.
\\\\
The failure of the global existence and decay property, which occurs for some solution $u$ with $|||u|||=A$ for any $A>A_c$, is expressed by $\|u\|_{\se_{T^*(u(0))}} = +\infty$, where we define for $T>0$
$$\|u\|_{\et} = \left( \sup_{t\in (0,T)}\|u(t)\|^2_{\hdhalf(\rt)} + \|D^{{3/2}}u\|^2_{L^2(\rt \times (0,T))}\right)^\frac{1}{2}$$
(so $\|u\|_{\se_{T^*(u(0))}} < +\infty$ whenever $|||u|||<A_c$ -- and therefore also, by standard embeddings, $u\in L^5(\rt \times (0,T))$ so such solutions are smooth by the ``Ladyzhenskaya-Prodi-Serrin condition", see e.g. \cite{ess}).
\\\\
In Theorems \ref{thm:a} and \ref{thm:b}, we establish the existence of a solution (``critical element") $u_c$ with initial datum $u_{0,c}$ and $T^*(u_{0,c}) < +\infty$ such that $$|||u_c||| = \sup_{t\in[0,T^*(u_{0,c}))}\|u_c(t)\|_{\dot{H}^\frac{1}{2}} = A_c \qquad \textrm{and} \qquad \|u_c\|_{\se_{T^*(u_{0,c})}} = +\infty $$
and, further, that for any such solution, there exist $x(t)\in \rt$, $0<\lambda(t) \to +\infty$ as $t \nearrow T^*(u_{0,c})$ and $s(t) \in [t,T^*(u_{0.c}))$ such that, for $u=u_c$ and $u_0 = u_{0,c}$,
\begin{equation}\label{eq:kay}
K:=\left\{\frac{1}{\lambda(t)}u\left(\frac{x - x(t)}{\lambda(t)}, s(t) \right),\ t\in [0,T^*(u_0))\right\}
\end{equation}
is pre-compact in $L^3(\rt)$.
\\\\
The pre-compactness in $\lt$ of such a $K$ is shown to be inconsistent with $T^*(u_0) < +\infty$ in Theorem \ref{thm:c}, by using the backwards uniqueness results for parabolic equations established in \cite{ess2,ess3} as well as the theory of unique continuation for parabolic equations to show that in fact $u_c \equiv 0$.  (In general, Step 3 requires something specific to the particular case being studied -- for example the Morawetz-type estimate for NLS used in \cite{km3} -- as opposed to the methods used to establish Steps 1 and 2 which are fairly general in nature.)
\\\\
One interesting scenario in which such a $K$ would be compact in $\lt$ is if, in (\ref{eq:kay}), one could take $\l(t) = (T^*-t)^{-\frac{1}{2}}$ for some $u_0 \in \lt$ (with $T^* =T^*(u_0) < +\infty$), $x(t) \equiv 0$ and $s(t) = t$, and one imposes that $K=\{U\}$ (i.e., $u(x,t) = \frac{1}{\sqrt{T^*-t}}U(\frac{x}{\sqrt{T^*-t}})$) for some given non-zero $U \in L^3$.  This is the case of a ``self-similar" solution which was first ruled out in the important paper \cite{necas} by J. Ne{\v c}as, M. R{\r u}{\v z}i{\v c}ka and V. {\v S}ver{\' a}k  in the (spatially global) $\lt$ setting (see \cite{tsai} for more general results).  Theorem \ref{thm:c} can therefore be thought of as a generalization of the result in \cite{necas}.  (Such solutions are of course ruled out as well by the more recent paper \cite{ess}, but the proof here is much simpler for that purpose.)
\\\\
\textbf{Acknowledgements:}  We would like to thank Professors Isabelle Gallagher and Vladim\' ir \v Sver\' ak for many helpful comments and discussions.  The second author would like to extend a personal and special thank you as well to Vladim\' ir \v Sver\' ak for his continuing help, guidance and support.
\section{Preliminaries}
We'll say that $u$ is a ``mild" solution of NSE on $[t_0,t_0 + T]$ for some $t_0 \in \R$ and $T>0$ if, for some divergence-free initial datum $u_0$, $u$ solves (in some function space) for $t\in [t_0,t_0+T]$ the integral equation
\begin{equation}\label{eq:mild}
u(t) = e^{(t-t_0)\Delta }u_0 + \int_{t_0}^t e^{(t-s)\Delta} \P\nabla \cdot (-u(s)\otimes u(s)) \ ds\ .
\end{equation}
We have used the following notation: For a tensor $F=(F_{ij})$ we define the vector  $\nabla \cdot F$ by $(\nabla \cdot F)_i = \sum_{j} \partial_j F_{ij}$, and for vectors $u$ and $v$, we define their tensor product $u\otimes v$ by $(u\otimes v)_{ij} = u_i v_j$.
\\\\
We'll consider solutions in spatial dimension three ($x\in \rt$), so $u=(u_1,u_2,u_3)$, $u_i = u_i(x,t)$, $1\leq i \leq 3$.  In that case, the projection operator $\P$ onto divergence-free fields is defined on a vector field $f$ by $(\P f)_j = f_j + \sum_{k=1}^3 \mathcal{R}_j
\mathcal{R}_kf_k$, $1\leq j \leq 3$, and the Riesz transform $\RR_j$ is defined on a scalar $g$ via Fourier transforms by $(\RR_j g)^{\wedge} (\xi) = \frac{i \xi_j}{|\xi|}\hat g (\xi)$.  (One can also formally write this as $\P f = f-\nabla \frac{1}{\Delta}(\nabla \cdot f)$.)  The heat kernel $e^{t\Delta}$ is defined by $e^{t\Delta}g(x) = [e^{-|\cdot|^2 t}\hat g (\cdot)]^\vee (x)= ((4\pi t)^{-3/2} \exp\{-|\cdot|^2/{4t} \} * g)(x)$, and extended to act component-wise on vector fields.
\\\\
Formally, (\ref{eq:mild}) comes from applying $\mathbb{P}$ to the classical Navier-Stokes equations which one can write as
\begin{equation}\label{eq:nse}
\begin{array}{rcl}
u_t -\D u + \nabla p & = &  \nabla \cdot (-u\otimes u) \\
\nabla \cdot u & = & 0
\end{array}
\end{equation}
(since $\nabla \cdot(u\otimes u) = (u\cdot \nabla )u$ due to the condition $\nabla \cdot u = 0$) and solving the resulting heat equation (since $\P (\nabla p) = 0$) by Duhamel's formula.
\\\\
In what follows, we'll set $L^p = L^p(\rt)$ and  $\|g\|_p = \|g\|_{L^p}$ for any $p\in [1,+\infty]$, $\dot H^s = \dot H^s (\rt) = \{g\in \mathscr{S}' \ | \ (D^s g)^{\wedge}(\cdot) = |\cdot|^\frac{s}{2}\hat{g}(\cdot) \in L^2\}$ for any $s\in \R$ where $\mathscr{S}'$ denotes the space of tempered distributions, and for $f=f(x,t)$ and any Banach space $X$, $f\in L^p((a,b);X) \iff \|f(t)\|_X = \|f(\cdot,t)\|_X\in L^p(a,b)$.  For any collection of Banach spaces $(X_m)_{m=1}^M$ and $X := X_1 \cap \cdots \cap X_M$, we'll always set $\|g\|_{X} = (\sum_{m=1}^M \|g\|^2_{X_m})^\frac{1}{2}$.  Similarly, for vector-valued $f=(f_1,\dots,f_M)$, we define $\|f\|_{X} = (\sum_{m=1}^M \|f_m\|^2_{X})^\frac{1}{2}$.  For easy reference, we collect and state here the definitions of the main spaces to which we refer throughout the paper:
$$
\mathscr{E}_T = L^\infty((0,T);\hdhalf) \cap L^2((0,T);\hdthalf)\ ;
$$
$$
\mathcal{E}^{(1/2)}_{T} :=\mathcal{C}([0,T);\hdhalf) \cap L^2 ((0,T);\hdthalf)\ ; \qquad
\mathcal{E}^{(3)}_{T}:=\mathcal{C}([0,T);L^3) \cap L^5 (\rt \times (0,T))\ .
$$
\section{Local theory}
In this section, we consolidate various known local existence results for the Navier-Stokes equations.  We also unify the various theories through known ``persistency" results.  Although the results presented in this section are well-known to experts, it seems to us that simple, self-contained proofs are often difficult to locate, so we present them here for the convenience of the reader.  Our main results are given in Section 3, and the expert reader may prefer to skip directly to that section now.
\\\\
The goal in what follows is to establish the existence of ``local" solutions to (\ref{eq:mild}) in some (space-time) Banach space $\X = \X_T$ of functions defined on $\rt \times [0,T)$ for some possibly small $T>0$, with divergence-free initial datum $u_0$ in a Banach space $X$.  In what follows, we will let $X$ equal $L^3$ or $\hdhalf$.  (See, e.g., \cite{cannone} and \cite{kt} respectively for local well-posedness in  $\dot B^{-1+\frac{3}{p}}_{p,\infty}$ and $BMO^{-1}$, and the recent ill-posedness result \cite{bp} for $\dot B^{-1}_{\infty,\infty}$.)  We'll re-write (\ref{eq:mild}) as
\begin{equation}\label{eq:abs}
x = y + B(x,x)
\end{equation}
and, under the assumption that $y\in \X$, try to solve the equation for some $x\in \X$ (where $\X$ will be chosen so that $u_0 \in X$ implies $\etl u_0 \in \X$).  This will be accomplished by the following abstract lemma, using the contraction mapping principle:
\begin{lemma}\label{lemma:abs}
Let $\X$ be a Banach space with norm $\|x\| = \|x\|_\X$, and let $B:\X \times \X \to \X$ be a continuous bilinear form such that there exists $\eta = \eta_\X >0$ so that
\begin{equation}\label{eq:f}
\|B(x,y)\| \leq \eta \|x\| \|y\|
\end{equation}
for all $x$ and $y$ in $\X$.  Then for any fixed $y\in \X$ such that $\|y\| < 1/{(4\eta)}$, equation (\ref{eq:abs}) has a
unique\footnote{In fact, the uniqueness can be improved to the larger ball of radius $\frac{1}{2\eta}$, see e.g. \cite{cannone2}, formula (122).} solution $\bar{x}\in \X$ satisfying $\|\bar{x}\| \leq R$, with
\begin{equation}\label{eq:g}
R:= \frac{1-\sqrt{1-4\eta \|y\|}}{2\eta} > 0\ .
\end{equation}
\end{lemma}
\noindent
{\bf Proof:} \quad Let $F(x) = y + B(x,x)\ .$  Using (\ref{eq:f}) and the triangle inequality, one can verify directly that $F$ maps $B_R := \{x\in \X \ | \ \|x\| \leq R\}$ into itself.  Moreover, $F$ is a contraction on $B_R$ as follows:
Suppose $x,x' \in B_R$.  Then
$$
\|F(x) - F(x')\|  =  \|B(x,x) - B(x',x')\|=  \|B(x-x',x) + B(x',x-x')\| \leq$$
$$ \leq  \eta \|x-x'\| \|x\|+ \eta \|x'\|\|x-x'\|
 \leq  2\eta R \|x-x'\| \ ,
$$
and clearly $2\eta R < 1$ by (\ref{eq:g}).
Hence the contraction mapping principle guarantees the existence of a unique fixed-point $\bar{x} \in B_R$ of the mapping $F$ satisfying $F(\bar{x}) = \bar{x}$, which proves the lemma. \hfill $\Box$
\subsection{Local theory in $\hdhalf$}
\footnote{This version of the local theory for initial data in $\hdhalf$ can be found in \cite{lr}.  For other versions, see for example the classical paper \cite{fk} and a more modern exposition in \cite{cannone}.      }Suppose $u_0 \in \hdhalf$, and let $\mathscr{E}_T = L^\infty((0,T);\hdhalf) \cap L^2((0,T);\hdthalf)$, with norm
\begin{equation}\label{eq:defet}
\|f\|_{\et} = \left( \|f\|^2_{L^\infty((0,T);\hdhalf)} + \|f\|^2_{L^2((0,T);\hdthalf)}\right)^\frac{1}{2}\ .
\end{equation}
Note that $\et \subset \ft :=L^4((0,T);\dot H^1)$, since H\" older's inequality gives
\begin{equation}\label{eq:h}
\|f\|_{\ft} \leq \|f\|^{\frac{1}{2}}_{L^\infty((0,T);\hdhalf)} \|f\|^\frac{1}{2}_{L^2((0,T);\hdthalf)} \ .
\end{equation}
By the well-known (see, e.g. \cite{chemin}, p. 120) inequality $\|h_1 h_2\|_{\hdhalf} \lesssim \|h_1\|_{\hdo}\|h_2\|_{\hdo}$, we have
$$\int_0^T \|f\otimes g\|^2_{\hdhalf}\ ds \lesssim \int_0^T \|f\|^2_{\hdo}\|g\|^2_{\hdo}\ ds \leq \left\{ \int_0^T \|f\|^4_{\hdo} \right\}^\frac{1}{2}
\left\{ \int_0^T \|g\|^4_{\hdo}\ ds \right\}^\frac{1}{2}\
$$
and hence
\begin{equation}\label{eq:i}
\|f\otimes g\|_{L^2((0,T);\hdhalf)} \lesssim \|f\|_{\ft}\|g\|_{\ft}\ .
\end{equation}
Denoting
$$B(f,g)(t):= \int_0^t e^{(t-s)\Delta} \P\nabla \cdot (-f(s)\otimes g(s)) \ ds\ ,$$
we can write
$$D^{\frac{3}{2}}B(f,g)(t) =  \int_0^t e^{(t-s)\Delta} \D F(s) \ ds\ $$
where $F(s):= \P\nabla \cdot (-\D)^{-1}D^{\frac{3}{2}}(f(s)\otimes g(s))$.  Now by the maximal regularity theorem for $\etl$ (see, e.g., \cite{lr}, Theorem 7.3), we have
$$\|D^{\frac{3}{2}}B(f,g)\|_{L^2(\rt \times (0,T))} \lesssim \|F\|_{L^2(\rt \times (0,T))}\ ,$$
and so since $F\sim D^{\frac{1}{2}}(f\otimes g)$, (\ref{eq:i}) gives
\begin{equation}\label{eq:j}
\|B(f,g)\|_{L^2((0,T);\hdthalf)} \lesssim \|f\|_{\ft}\|g\|_{\ft}\ .
\end{equation}
Let's recall the following lemma (see \cite{lr}, Lemma 14.1):
\begin{lemma}\label{lemma:lr}
let $T \in (0,+\infty]$ and $1\leq j \leq 3$.  If $h\in L^2(\rt \times (0,T))$, then $\int_0^t e^{(t-s)\D}\partial_j h \ ds \in \mathcal{C}_b([0,T);L^2)$.
\end{lemma}
\noindent
($\mathcal{C}_b$ indicates bounded continuous functions.)  For $f,g\in \ft$, (\ref{eq:i}) shows that $D^\frac{1}{2} (f\otimes g) \in L^2(\rt \times (0,T))$, so Lemma \ref{lemma:lr} gives $D^\frac{1}{2} B(f,g) \in \mathcal{C}([0,T);L^2)$ and hence $B(f,g) \in \mathcal{C}([0,T);\hdhalf)$.  Moreover, the proof of Lemma \ref{lemma:lr} also gives the estimate
$$\|B(f,g)\|_{L^\infty((0,T);\hdhalf)} \lesssim \|f\otimes g\|_{L^2((0,T);\hdhalf)} \ ,$$ which, together with (\ref{eq:i}) gives
\begin{equation}\label{eq:k}
\|B(f,g)\|_{L^\infty((0,T);\hdhalf)} \lesssim  \|f\|_{\ft}\|g\|_{\ft}\ .
\end{equation}
Using (\ref{eq:h}), (\ref{eq:j}) and (\ref{eq:k}), we now conclude that
\begin{equation}\label{eq:l}
\|B(f,g)\|_{\ft} \leq \eta  \|f\|_{\ft}\|g\|_{\ft}\
\end{equation}
for some $\eta >0$.  (We remark that $\eta$ is independent of $T$.)  By the standard $L^2$ energy estimates for the heat equation, $u_0 \in \hdhalf$ implies that $U \in \et \subset \ft$, where $U(t):= \etl u_0$.  Therefore we can take $T$ small enough that $\|U\|_{\ft} < \frac{1}{4\eta}$ (or, if $\|u_0\|_\hdhalf$ is small enough one may take $T=+\infty$, giving a ``small data" global existence result), and hence by Lemma \ref{lemma:abs} with $\X = \ft$, there exists a unique small mild solution $u\in \ft$ of NSE on $[0,T)$.  Note that (\ref{eq:j}) and (\ref{eq:k}) show that $B(u,u)\in \et$ as well, so that more specifically $u = U + B(u,u) \in \et$.  Moreover, $u\in \mathcal{C}([0,T);\hdhalf)$ by Lemma \ref{lemma:lr} and the standard theory for the heat equation (see, e.g.,  \cite{ural}).
\subsection{Local theory in $\hdhalf \cap L^\infty$}
The above result can be refined to give a solution which not only remains in $\hdhalf$, but belongs to $L^\infty$ as well for $t>0$.  This will be shown by considering the spaces
$$\eti := \et \cap \{ g\ |\ \sqrt{t}g(x,t) \in L^\infty(\rt \times (0,T))\}$$
$$\fti := \ft \cap \{ g\ |\ \sqrt{t}g(x,t) \in L^\infty(\rt \times (0,T))\}\ .$$
We will show that there exists some $\eta_{\infty} >0$ such that
\begin{equation}\label{eq:m}
\|B(f,g)\|_{\fti} \leq \eta_{\infty} \|f\|_{\fti}\|g\|_{\fti}\ ,
\end{equation}
and hence there exists a unique small solution $u\in \fti$ by Lemma \ref{lemma:abs} so long as $U(t) = \etl u_0$ satisfies
\begin{equation}\label{eq:n}
\|U\|_{\fti} < \frac{1}{4\eta_\infty}\ .
\end{equation}
Note that Young's inequality gives
\begin{equation}\label{eq:o}
\|e^{t\Delta}u_0\|_\infty \lesssim t^{-1/2} \|u_0\|_3\ ,
\end{equation}
hence $u_0 \in \hdhalf \subset L^3$ implies that indeed $U\in \fti$.  Moreover, as before, the resulting solution in the case of (\ref{eq:n}) will belong more specifically to $\eti$.
\\\\
We claim now that
\begin{equation}\label{eq:p}
\lim_{t\to 0} \|\sqrt{t}\etl u_0\|_\infty = 0 \quad \forall\ \  u_0 \in L^3\ ,
\end{equation}
which will show that (\ref{eq:n}) will hold for any $u_0 \in \hdhalf$ for $T$ sufficiently small.  To prove (\ref{eq:p}), for any $\e >0$, take $R$ and $M$ large enough that $\|u_0(1-\chi^{M,R})\|_3 < \e/2$,
where $\chi^{M,R}(x) = 1$ for $x\in \{ |x|<R \} \cap \{x \ | \ |u_0(x)| \leq M\}$ and $0$ otherwise.  Then by Young's inequality we have
$$\|\sqrt{t}\etl u_0\|_\infty \leq \|\sqrt{t}\etl (u_0\chi^{M,R})\|_\infty + \|\sqrt{t}\etl (u_0(1-\chi^{M,R}))\|_\infty \leq \sqrt{t}\|\etl\|_1 \cdot M + \frac{\e}{2} < \e$$
for small enough $t>0$.  The bilinear estimate (\ref{eq:m}) is a consequence of the continuous embedding $\hdhalf \subset L^3$, estimates (\ref{eq:k}) and (\ref{eq:l}) and the following claim:
\begin{clm}\label{clm:b}
$$\sup_{t\in(0,T)} \sqrt{t}\|B(f,g)(t)\|_\infty \lesssim \|B(f,g)\|_{L^\infty((0,T);L^3)} + \sup_{t\in(0,T)}\sqrt{t} \|f(t)\|_\infty \cdot \sup_{t\in(0,T)} \sqrt{t}\|g(t)\|_\infty$$
\end{clm}
\noindent
{\bf Proof:} \quad We'll need the following facts:
\begin{equation}\label{AA}
\begin{array}{rl}
(i) & \|\etl u_0\|_\infty \lesssim t^{-\frac{1}{2}} \|u_0\|_3 \\\\
(ii) & \etl \P \nabla \cdot h = \frac{1}{t^2}H\left(\frac{\cdot}{\sqrt{t}} \right) * h,\ \textrm{where}\  |H(y)| \lesssim (1+|y|)^{-4}.
\end{array}
\end{equation}
(i) is just Young's inequality, and (ii) can be found for example in \cite{lr}, Prop. 11.1 on ``The Oseen Kernel".  (See also \cite{solonnikov}, translated in \cite{soltr}).  Now write
$$B(f,g) = e^{(t/2)\D}B(f,g)(t/2) + \int_{t/2}^t e^{(t-s)\D}\P\nabla \cdot (f\otimes g)(s)\ ds\ .$$
By (ii) and a change of variables we have
$$|e^{(t-s)\D}\P \nabla \cdot (f\otimes g)(x,s)| = \left|\frac{1}{(t-s)^{1/2}} \int_{\rt} H(z)(f\otimes g)(x-z\sqrt{t-s},s)\ dz \right|$$
$$\qquad \qquad \lesssim \frac{\|f(s)\|_\infty \|g(s)\|_\infty}{(t-s)^{1/2}} \int_{\rt} \frac{1}{(1+|z|)^4}\ dz  = C \cdot \frac{\|f(s)\|_\infty \|g(s)\|_\infty}{(t-s)^{1/2}}\ ,$$
and hence by (i)
$$
\begin{array}{rcl}
\|B(f,g)(t)\|_\infty & \lesssim & (t/2)^{-1/2} \|B(f,g)(t/2)\|_3 + \int_{t/2}^t (t-s)^{-1/2}\|f(s)\|_\infty \|g(s)\|_\infty \ ds        \\\\
& \lesssim & \displaystyle{t^{-1/2} \|B(f,g)(t/2)\|_3 +  t^{1/2} \sup_{s\in (\frac{t}{2},t)}\|f(s)\|_\infty \cdot \sup_{s\in (\frac{t}{2},t)}\|g(s)\|_\infty \ .}
\end{array}
$$
Therefore, for $t\in (0,T)$,
$$
t^{1/2}\|B(f,g)(t)\|_\infty  \lesssim   \|B(f,g)(t/2)\|_3 +  t^{1/2} \sup_{s\in (\frac{t}{2},t)}\|f(s)\|_\infty \cdot t^{1/2}\sup_{s\in (\frac{t}{2},t)}\|g(s)\|_\infty  $$$$
 \leq  \sup_{t\in(0,T)} \|B(f,g)(t)\|_3 + 2 \sup_{t\in(0,T)}t^{1/2} \|f(t)\|_\infty \cdot \sup_{t\in(0,T)}t^{1/2}\|g(t)\|_\infty
$$
which proves the claim.  (We remark that the constant in the claim does not depend on $T$.)
\subsection{Local theory in $\lt$ (and $L^3 \cap L^\infty$)}
Take $u_0\in \lt$ and let $D_T = \rt \times (0,T)$ for some $T\in (0,+\infty]$.  We will show local existence\footnote{This version of the local theory for initial data in $L^3$ was presented in a course on mathematical fluid mechanics given by Prof. Vladim\' ir \v Sver\' ak at the University of Minnesota in the spring of 2006, and can be found as well in \cite{dongdu}.  For other versions, see e.g. the classical paper \cite{kato} or the more modern treatment in \cite{cannone}.         } in the space $L^5(D_T)$.  Note first that $U\in L^5(D_T)$ where $U(x,t) = \etl u_0(x)$ as follows:
Since $\etl u_0 = \etl (u_0)^+ - \etl (u_0)^-$\ , we can assume $u_0 \geq 0$.  Since $$(U_t - \D U)\cdot U^{p-1} = 0\ $$
for any $p>1$, integration by parts yields
$$
\int_{\rt} |U^{p/2}(x,t)|^2\ dx + \frac{4(p-1)}{p}\int_0^t\int_{\rt} |\nabla(U^{p/2})(x,t)|^2\ dx\ dt = \int_{\rt} |U^{p/2}(x,0)|^2\ dx \ .
$$
Taking $p=3$ and using the inequality
$$\|g\|_{L^{10/3}(D_T)} \lesssim \|g\|^{2/5}_{L^\infty((0,T);L^2)} \|g\|^{3/5}_{L^2((0,T);\hdo)}$$
which is due to the H\" older\footnote{Note that H\" older's inequality implies the following interpolation inequality:  for $p< r < q$, $\alpha = \frac{\frac{1}{r} - \frac{1}{q}}{\frac{1}{p} - \frac{1}{q}}\in (0,1)$ and $\|g\|_r = \||g|^\alpha |g|^{1-\alpha}\|_r \leq \|g\|^\alpha_p \|g\|^{1-\alpha}_q$.       } and Sobolev inequalities, we have
$$\|U^{3/2}\|_{L^{10/3}(D_T)} \lesssim \|U^{3/2}(0)\|_2$$
which is exactly
\begin{equation}\label{eq:five}
\|U\|_{L^5(D_T)} \lesssim \|u_0\|_3\ .
\end{equation}
Recall now that
$$B(f,g)(x,t) = \int_0^t K(\cdot, t-s) * (f(\cdot, s)\otimes g(\cdot, s))(x)\ ds\ ,$$
where we can write (see (\ref{AA}))
$$K(x,t) =
\left\{
\begin{array}{ccc}
\frac{1}{t^2}H\left(\frac{x}{\sqrt{t}}\right) & , & t>0 \\
0 & , & t\leq 0
\end{array}
\right.
$$
for some smooth $H\in L^1 \cap L^\infty$.  With a slight abuse of notation for simplicity, we now make the following claim:
\begin{clm}\label{clm:c}
$$\left\|\int_0^t K(t-s) * h(s)\ ds \right\|_{L^5(\rt \times \R)} \lesssim \|h\|_{L^{5/2}(\rt \times \R)}\ ,$$  whenever the right-hand side is finite.
\end{clm}
\noindent
{\bf Proof:}
$$\left\| \left\|\int_0^t K(t-s)*h(s)\ ds\right\|_{L^5_x} \right\|_{L^5_t(\R)}\ \
\leq \ \ \left\| \int_0^t \left\|K(t-s)*h(s)\right\|_{L^5_x}\ ds \right\|_{L^5_t(\R)} $$
$$
\begin{array}{rcl}
& \leq  & \displaystyle{ \left\| \int_0^t \|K(t-s)\|_{L^{5/4}_x}\|h(s)\|_{L^{5/2}_x}\ ds \right\|_{L^5_t(\R)}  } \\\\
& =  & \displaystyle{ \left\| \int_0^t (t-s)^{-4/5} \|H\|_{L^{5/4}_x}\|h(s)\|_{L^{5/2}_x}\ ds \right\|_{L^5_t(0,+\infty)}  } \\\\
& \lesssim  & \displaystyle{ \left\| \int_{-\infty}^{+\infty}  \frac{\|h(s)\|_{L^{5/2}_x}}{|t-s|^{1-\alpha}}\ ds \right\|_{L^5_t(\R)}  \ \ \leq \ \
 \left\| \|h(s)\|_{L^{5/2}_x} \right\|_{L^{5/2}_t(\R)}  } \\\\
\end{array}
$$
where $\alpha = 1/5$ and we have used Young's inequality\footnote{$\|f*g\|_r \leq \|f\|_p\|g\|_q$ for $1\leq p,q,r \leq +\infty$ whenever $\frac{1}{p}+\frac{1}{q}-1=\frac{1}{r}>0$} in the $x$ variable with
$\frac{1}{5} = \frac{1}{(5/4)} + \frac{1}{(5/2)} - 1$
in the second line, and one-dimensional fractional integration\footnote{In dimension $n$, $\|D^{-\alpha}f\|_q = C_\alpha \|\ |\cdot|^{-n+\alpha}*f(\cdot)\ \|_q \leq C_{p,q}\|f\|_q$ for $1 < p < q < +\infty$ and $0< \alpha < n$ whenever $\frac{1}{q} = \frac{1}{p} - \frac{\alpha}{n}$ (see \cite{stein2})}  in the $t$ variable with $\frac{1}{5} = \frac{1}{(5/2)} - \alpha$ in the last.
Now since for $t\in (0,T)$ we have
$$\int_0^t K(t-s)*h(s)\ ds = \int_0^t K(t-s)* (f\otimes g)(s)\ ds$$
for $h:= \chi_{[0,T]} f\otimes g$ where $\chi_{[0,T]}$ is the indicator function for the interval $[0,T]$, for any fields $f,g \in L^5(D_T)$ Claim \ref{clm:c} gives
\begin{equation}\label{eq:q}
\|B(f,g)\|_{L^5(D_T)} \leq \eta_5 \|f\otimes g\|_{L^{5/2}(D_T)} \leq \eta_5 \|f\|_{L^5(D_T)}\|g\|_{L^5(D_T)}
\end{equation}
for some $\eta_5 > 0$ independent of $T$.  Hence there exists a unique solution $u\in L^5(D_T)$ by Lemma \ref{lemma:abs} so long as $\|U\|_{L^5(D_T)} < \frac{1}{4\eta_5}$ which is true for small enough $T$ due to (\ref{eq:five}) (or, for $\|u_0\|_3$ sufficiently small, one may take $T=+\infty$).
Moreover, estimate (7.26) of \cite{ess} implies\footnote{Strictly speaking, the estimate was derived under the assumption that $u$ satisfies the differential NSE for a suitable pressure $p$; however since the mild solution is smooth, we can reduce to this case, see the proof of Theorem \ref{thm:c}.} that
\begin{equation}\label{eq:r}
\|B(f,g)\|_{L^\infty((0,T);L^3)} \lesssim \|f\|_{L^5(D_T)}\|g\|_{L^5(D_T)}\ ,
\end{equation}
the proof of which shows moreover that $t \longmapsto \|B(f,g)(t)\|_3$ is continuous.  This, along with the fact that $B(f,g)(t)$ is weakly continuous in $\lt$ on $(0,T)$ (see \cite{ess}, (7.27)) implies that $B(f,g)\in \mathcal{C}((0,T);L^3)$.  The standard theory of the heat equation implies now that $u = U + B(u,u) \in \mathcal{C}([0,T); L^3)$.
\begin{remark}\label{remark:bdd}
We can furthermore assume that $\sqrt{t}u(x,t) \in L^\infty(D_T)$ (for a possibly smaller $T$) in a manner similar to the proof of (\ref{eq:m}), using Claim \ref{clm:b} and estimate (\ref{eq:r}).  For a construction of local solutions with higher regularity in similar spaces, see \cite{dongdu}.
\end{remark}
\noindent
At this point, we also mention the following uniqueness theorem for solutions in the class $\mathcal{C}([0,T);\lt)$ established in \cite{flt} (see also \cite{monniaux} for a simplified proof):
\begin{thm}\label{unique}
Let $u^1, u^2 \in \mathcal{C}([0,T);\lt)$ both satisfy (\ref{eq:mild}) for a fixed $u_0\in \lt$.  Then $u^1 \equiv u^2$ on $[0,T)$.
\end{thm}
\noindent
Note that there is no size restriction on $\|u_0\|_3$.
\subsection{Unification of the theories and further properties}
For any $u_0\in \lt$, the local theory guarantees the existence of a mild solution $u_{(T)} \in \mathcal{E}^{(3)}_{T}$ on $[0,T)$ for some $T>0$ with $u_{(T)}(0) = u_0$, where we define
\begin{equation}\label{ett}
\mathcal{E}^{(3)}_{T}:=\mathcal{C}([0,T);L^3) \cap L^5 (\rt \times (0,T)) \ .
\end{equation}
By Theorem \ref{unique}, there can be at most one mild solution in $\mathcal{C}([0,T);L^3)$ with initial datum $u_0$ for a fixed $T>0$, and hence there can be at most one\footnote{In fact, since $u\in \mathcal{E}^{(3)}_{T}$ contains the additional information that $u\in L^5_{x,t}$, one can derive uniqueness from the local theory without using Theorem \ref{unique}, but we proceed in this way for simplicity.} such $u_{(T)}$.  Define $T_{3}^* = T_3^*(u_0) > 0$ by
$$T_{3}^*(u_0) := \sup \{T>0 \ | \ \exists! \ \textrm{mild solution}\ u_{(T)}\in \mathcal{E}^{(3)}_{T}\ \textrm{on}\ [0,T)\ \mathrm{s.t.}\ u_{(T)}(0) = u_0 \}\ .$$
Again, by the $L^3$ uniqueness theorem, $u_{(T_1)}(0) = u_{(T_2)}(0)$ implies that $u_{(T_1)}(t) = u_{(T_2)}(t)$ for all  $t\in [0,\min\{T_1,T_2\})$, hence there exists a unique mild solution $u^{(3)}$ on $[0,T^*_3)$ such that $u^{(3)}(0) = u_0$ and $u^{(3)}\in \mathcal{E}^{(3)}_{T}$ for any $T<T^*_3$.  We'll call $T^*_3$ the maximal time of existence in $L^3$ of the mild solution associated to $u_0$.
\\\\
Similarly, for $u_0 \in \hdhalf$, defining
\begin{equation}\label{B}
\mathcal{E}^{(1/2)}_{T} :=\mathcal{C}([0,T);\hdhalf) \cap L^2 ((0,T);\hdthalf) \
\end{equation}
and $T_{1/2}^* = T_{1/2}^*(u_0) >0$ by
$$T_{1/2}^*(u_0):= \sup \{T>0 \ | \ \exists! \ \textrm{mild solution}\ u_{(T)}\in \mathcal{E}^{(1/2)}_{T}\ \textrm{on}\ [0,T)\ \mathrm{s.t.}\ u_{(T)}(0) = u_0 \}\ ,$$
there exists a unique mild solution $u^{(1/2)}$ on $[0,T^*_{1/2})$ such that $u^{(1/2)}(0) = u_0$ and $u^{(1/2)}\in \mathcal{E}^{(1/2)}_{T}$ for any $T<T^*_{1/2}$.  We'll call $T^*_{1/2}$ the maximal time of existence in $\hdhalf$ of the solution associated to $u_0$.
\\\\
The standard embeddings and interpolation inequalities give the estimate
$$\|g\|_{L^5(\rt \times (0,T))} \lesssim \|g\|^{3/5}_{L^\infty((0,T);\hdhalf)} \|g\|^{2/5}_{L^2((0,T);\hdthalf)}\ ,$$
and hence $\mathcal{E}^{(1/2)}_{T} \subset \mathcal{E}^{(3)}_{T}$ for any $T>0$.  This implies that $T^*_{1/2} \leq T^*_3$ for any $u_0 \in \hdhalf \subset L^3$.  We will now proceed to show that in fact $T^*_{1/2} = T^*_3$ for $u_0 \in \hdhalf$.
\\\\
The following lemma will show that a solution which is continuous in $\hdhalf$ for some time cannot leave $\hdhalf$ while remaining continuous in $L^3$ (or, equivalently, while remaining in $L^5_{x,t}$), and hence that the strict inequality $T^*_{1/2} < T^*_3$ is impossible for $u_0\in \hdhalf$.  Furthermore, Theorem \ref{unique} implies that $u^{(1/2)} = u^{(3)}$ on $[0,T^*)$ where $T^* = T^*_{1/2} = T^*_3$.  In what follows, for $u_0\in \lt$ (in particular, for $u_0 \in \hdhalf$) we'll denote $u^{(3)}$ ($=u^{({1/2})}$) $=:NS(u_0)$.
\begin{lemma}\label{lemma:a}
Suppose $u\in \mathcal{C}([0,T];\lt)$ is a mild solution for NSE on $[0,T]$.  Then
\begin{itemize}
\item[(a)] $u\in L^5 (\rt \times (0,T))$;
\item[(b)] If moreover $u(0)=u_0 \in \hdhalf$, then $u\in \mathcal{C}([0,T];\hdhalf) \cap L^2((0,T);\hdthalf)$.
\end{itemize}
\end{lemma}
\noindent
In order to prove the lemma, we will need the following claim:
\begin{clm}\label{clm:a}
Suppose $u\in \mathcal{C}([0,T];\lt)$, and fix some $a>0$.  Then
\begin{itemize}
\item[(i)] $s\longmapsto \sup_{t\in (0,a)} \sqrt{t}\|e^{t\Delta}u(s)\|_\infty \in \mathcal{C}[0,T]$\  and
\item[(ii)] $s\longmapsto \|e^{t\Delta}u(x,s)\|_{L^5_{x,t}(\rt \times (0,a))} \in \mathcal{C}[0,T]$\ .
\end{itemize}
If, moreover, $u\in \mathcal{C}([0,T];\hdhalf)$, then
\begin{itemize}
\item[(iii)] $s\longmapsto \|e^{t\Delta}u(s)\|_{L^4_t((0,a);\dot H^1)} \in \mathcal{C}[0,T]$\ .
\end{itemize}
\end{clm}
\ \\
{\bf Proof of Claim \ref{clm:a}:}
\quad
The claim is a simple consequence of the linearity of the operator $\etl$, estimates (\ref{AA}(i)) and (\ref{eq:five}), and the estimate $\|e^{t\Delta}g\|_{L^4((0,a);\dot H^1)} \lesssim \|g\|_{\hdhalf}$ which follows from the $L^2$ energy inequality for the heat equation and the standard embeddings.  Writing any one of these as $\|\etl g\|_{\X_a} \lesssim \|g\|_X$ where $u\in \mathcal{C}([0,T];X)$, we have
$$|\ \|\etl u(s_2)\|_{\X_a} - \|\etl u(s_1)\|_{\X_a}\ | \leq \|\etl (u(s_2)-u(s_1))\|_{\X_a} \lesssim \|u(s_2) - u(s_1)\|_X$$
which proves the claim. \hfill $\Box$
\\\\
{\bf Proof of Lemma \ref{lemma:a}:}
\quad
Note first that we may actually assume $u$ is a solution in $\mathcal{C}([0,T+\epsilon);L^3)$ on $[0,T+\epsilon)$ for some $\epsilon >0$.  This is guaranteed by the local theory in $L^5_{x,t}$ for initial data (e.g., $u(T)$) in $L^3$.  By the same local existence theory, we are guaranteed a solution $u^{(3)} \in \mathcal{E}^{(3)}_{T_1} =\mathcal{C}([0,T_1);L^3) \cap L^5 (\rt \times (0,T_1))$ with initial datum $u(0)$ for some small $T_1>0$.  By the uniqueness guaranteed by Theorem \ref{unique}, we know that in fact $u^{(3)} = u$ on $[0,T_1)$.
\\\\
Set $\overline{T}:= \sup \{T'\in (0,T+\e)\ |\ u\in \mathcal{E}^{(3)}_{T'}\}$.  If $\overline{T} > T$ then we have proved (a), so suppose $\overline{T} \leq T$.
Since $u(\TT) \in \lt$, (\ref{eq:five}) shows that $e^{t\Delta} u(\TT)\in L^5(\rt \times (0,+\infty))$, hence one may take some $t_1 \in (0,\e)$ small enough so that
\begin{equation}\label{eq:a}
\|e^{t\Delta} u(\TT)\|_{L^5(\rt \times (0,t_1))} < \frac{\e_0}{2}\ ,
\end{equation}
where $\e_0 >0$ is such that $\|e^{t\Delta} v\|_{L^5(\rt \times (0,t_1))} < \e_0$ guarantees an $\mathcal{E}^{(3)}_{t_1}$-solution on $(0,t_1)$ with initial datum $v$ by Lemma \ref{lemma:abs} and (\ref{eq:q}). Now, by Claim \ref{clm:a} (ii), there exists $\delta \in (0,t_1)$ such that
\begin{equation}\label{eq:b}
\|e^{t\Delta} u(\TT-\delta)\|_{L^5(\rt \times (0, t_1 ))} < \e_0\ .
\end{equation}
Therefore by the local existence theory, we are guaranteed a solution $u^{(3)} \in \mathcal{C}([\TT- \delta, \TT + \overline{\delta});L^3)\cap L^5(\rt \times (\TT - \delta,\TT + \overline{\delta}))$, where $\overline{\delta}:= t_1 - \delta >0$.  Again by Theorem \ref{unique}, $u^{(3)}=u$ on $(\TT- \delta, \TT + \overline{\delta})$.  But now clearly $u\in L^5 (\rt \times (0,\TT + \overline{\delta}))$, which is a contradiction to the definition of $\TT$, and (a) is proved.
\\\\
Let's now turn to (b).  First suppose that we can show that $u_0 \in \hdhalf$ implies $u\in \mathcal{C}([0,T];\hdhalf)$. Then by replacing $\mathcal{E}^{(3)}_T$ above by $\mathcal{F}^{( 1/2)}_T := \mathcal{C}([0,T);\hdhalf) \cap L^4((0,T);\dot H^1)$ and replacing Claim \ref{clm:a} (ii) by Claim \ref{clm:a} (iii),
the method used to prove (a) can be adjusted to show that $u\in L^4((0,T);\dot H^1)$. This fact, along with standard estimates (see the proof of the local theory in $\hdhalf$), yields the conclusion in (b).
It therefore remains only to show that
\begin{equation}\label{eq:c}
u\in \mathcal{C}([0,T];\lt), u_0 \in \hdhalf \ \   \Longrightarrow \ \ u\in \mathcal{C}([0,T];\hdhalf)\ .
\end{equation}
In order to do this, we turn to the local solutions $u^{(KT)}$ of Koch-Tataru \cite{kt} for (possibly large, but somewhat regular) initial data in the space $bmo^{-1}$, and a ``persistency" theorem given in \cite{flzz} which states that whenever the initial datum belongs to various spaces of regular data, including $\hdhalf$, the solution with that datum remains there for positive times.\footnote{We remark that it should be possible to complete the proof within the framework of the $\lt$-theory, without resorting to the stronger tools in \cite{kt} and \cite{flzz} which we have used here only in the interest of expediency.  }  To complete the proof, we rely as before on the uniqueness provided by Theorem \ref{unique}.
\\\\
We recall briefly the space $\mathcal{E}^{(KT)}_T$ in which the Koch-Tataru solutions were constructed:
$$\mathcal{E}^{(KT)}_T:=\{f\ |\ \sup_{t\in (0,T)} \sqrt{t}\|f(t)\|_\infty < +\infty\}\cap \qquad \qquad \qquad \qquad  \qquad \qquad  \qquad \qquad$$
$$ \qquad \qquad \qquad \qquad  \qquad \cap\{f\ |\ \sup_{x_0\in \rt, t\in (0,T)} \frac{1}{t^{3/2}} \int_0^t \int_{B_{\sqrt{t}}(x_0)} |f(x,s)|^2\ dx\ ds < +\infty\}\ ,$$
where for $r>0$ we denote $B_r(x_0) = \{x\ | \ |x-x_0| < r\}$.  Recall also the definition (see, e.g., \cite{lr}) that $v\in bmo^{-1}$ if and only if
\footnote{Strictly speaking, one should denote this space by $bmo^{-1}_T$, but we consider some fixed $T>0$.}
$$\|v\|_{bmo^{-1}} := \sup_{x_0\in \rt, t\in (0,T)} \frac{1}{t^{3/2}} \int_0^t \int_{B_{\sqrt{t}}(x_0)} |e^{s\Delta}v(x)|^2\ dx\ ds\  < \ +\infty\ .$$
By (\ref{eq:o}), $e^{t\Delta}u_0 \in \mathcal{E}^{(KT)}_T$ for any $u_0 \in L^3 \subset bmo^{-1}$.  We know from (\ref{eq:p}) that $\lim_{t\to 0} \sqrt{t}\|e^{t\D}u_0\|_\infty =0$, and it can moreover be shown that $\lim_{\l \searrow 0} \|u_0^\l\|_{bmo^{-1}} =0$, where $v^\l(x) = \l v(\l x)$, for any $u_0 \in L^3$. Note also that $e^{t\D}v^\l(x) = \l (e^{\l^2t\D}v)(\l x)$, and hence by a change of variables we have
$$\sup_{t\in(0,T)} \sqrt{t}\|e^{t\D}u_0^\lambda\|_\infty = \sup_{t\in(0,\l^2T)} \sqrt{t}\|e^{t\D}u_0\|_\infty$$
which therefore tends to zero as well as $\l \searrow 0$ by the preceding statement.
\\\\
The local result in \cite{kt} is that for any $T>0$ there exists some small $\e_T >0$ such that $\|e^{t\D}v_0\|_{\mathcal{E}^{(KT)}_T} < \e_T$ implies the existence of a unique small solution $v^{(KT)} \in \mathcal{E}^{(KT)}_T$ on $[0,T)$ satisfying $v(0) = v_0$. By the statements above, for any $u_0 \in L^3$ one can take $\l_1 >0$ small enough so that $\sup_{t\in (0,1)}\sqrt{t}\|e^{t\D}u_0^{\l_1}\|_\infty$, $\|u_0^{\l_1}\|_{bmo^{-1}} < \e_1$.  Letting $v_0 = u_0^{\l_1}$ we obtain a solution $v$ on $[0,1)$ with initial datum $u_0^{\l_1}$.  Due to the natural scaling of the equation, $u^{(KT)}(x,t) := \frac{1}{\l_1} v(\frac{x}{\l_1}, \frac{t}{\l_1^2})$ is a solution on $[0,(\l_1)^2)$ satisfying $u^{(KT)}(0) = u_0$.
\\\\
We return now to the proof of (\ref{eq:c}).  Suppose $u_0\in \hdhalf$.  By the local existence theory, there exists a solution $u^{(1/2)} \in \mathcal{E}^{( 1/2)}_{T_1}$ on $[0,T_1)$ with $u^{(1/2)}(0) = u_0$ for some small $T_1 > 0$.  By uniqueness (Theorem \ref{unique}), $u^{(1/2)} = u$ on $[0,T_1)$ hence, in particular, $u\in \mathcal{C}([0,T_1);\hdhalf)$.
\\\\
Define $\TT_1 := \max \{\ T_1 \in (0,T+\e) \ | \ u\in \mathcal{C}([0,T_1);\hdhalf)\}$, and suppose $\TT_1 \leq T$.
\\\\
Since $u(\To) \in L^3 \subset bmo^{-1}$, there exists some $\l_1 > 0$ small enough that
$$\sup_{t\in (0,1)}\sqrt{t}\|e^{t\D}[u(\To)]^{\l_1}\|_\infty,\ \|[u(\To)]^{\l_1}\|_{bmo^{-1}} < \e_1/2\ ,$$
guaranteeing a solution $\ukt$ on $[\To, \To + t_1)$ satisfying $\ukt(\To) = u(\To)$, where $t_1 = (\l_1)^2$.  Due to the embedding $\|u\|_{bmo^{-1}} \lesssim \|u\|_3$, we have $u\in \mathcal{C}([0,T);bmo^{-1})$, and in particular there exists some $\ttil < \To$ such that for  $\d: = \To - \ttil$ small enough we have
$$
\begin{array}{rcl}
\|[u(\ttil)]^\l - [u(\To)]^\l\|_{bmo^{-1}} &\lesssim &\|[u(\ttil)]^\l - [u(\To)]^\l\|_{3} \\\\
& =& \|u(\ttil) - u(\To)\|_{3} \ < \ \e_1/2 \ ,
\end{array}$$
so that $\|[u(\ttil)]^\l\|_{bmo^{-1}} < \e_1$.    Also, for $\d $ small enough, we have by (i) of Claim \ref{clm:a}
$$
\sup_{t\in (0,1)}\sqrt{t}\|e^{t\D}[u(\ttil)]^{\l_1}\|_\infty =\sup_{t\in (0,t_1)}\sqrt{t}\|e^{t\D}u(\ttil)\|_\infty \leq
\qquad \qquad \qquad \qquad \qquad \qquad$$
$$
\qquad \stackrel{(i)}{\leq}  \sup_{t\in (0,t_1)}\sqrt{t}\|e^{t\D}u(\To)\|_\infty + \e_1/2 \\\\
=  \sup_{t\in (0,1)}\sqrt{t}\|e^{t\D}[u(\To)]^{\l_1}\|_\infty + \e_1/2
<  \e_1 .
$$
Hence we are also guaranteed a solution $\ukt$ on $[\ttil, \ttil + t_1)$ satisfying $\ukt(\ttil) = u(\ttil)$, where $t_1 = (\l_1)^2$.  Taking $\d < t_1$, we now have a Koch-Tataru solution $\ukt$ on $[\ttil,\To + \overline{\d})$ satisfying $\ukt(\ttil) = u(\ttil)$, where $\overline{\d} = t_1 - \d >0$.
\\\\
Since $\ttil < \To$, we have $u(\ttil) \in \hdhalf$ by our definition of $\To$.  Now, the persistency principle of \cite{flzz} states that a Koch-Tataru solution obtained by the fixed-point argument on an interval with initial datum in $\dot B^{s,q}_p$ for $1\leq p,q \leq +\infty,\ s>-1$ will remain in $\dot B^{s,q}_p$ on that interval and is moreover continuous in that norm\footnote{The theorem in \cite{flzz} includes persistence in inhomogeneous spaces and Lebesgue spaces as well, including $L^\infty$.}.  Since\footnote{These are the standard Besov and Triebel-Lizorkin spaces; see, e.g., \cite{grafakos}.} $\dot B^{\frac{1}{2},2}_2 = \dot F^{\frac{1}{2},2}_2 = \hdhalf$, we see that $\ukt \in \mathcal{C}([\ttil,\To + \overline{\d});\hdhalf)$.  Theorem \ref{unique} then implies that $\ukt = u$ on $[\ttil, \To + \overline{\d})$, and hence $u\in \mathcal{C}([0,\To + \overline{\d});\hdhalf)$.  This, however, contradicts our assumption on $\To$, hence necessarily $\To > T$, and (\ref{eq:c}) is proved, completing the proof of Lemma \ref{lemma:a}. \hfill $\Box$
\\\\
Finally, we state here for convenience some known a priori results regarding globally-defined mild solutions to NSE evolving from possibly large data:
\begin{thm}[Decay]\label{decay}
Let $u\in \mathcal{C}([0,\infty);X)$ be a global-in-time solution to (\ref{eq:mild}) for some divergence-free $u_0\in X$, where $X$ is either $\hdhalf(\rt)$ or $L^3(\rt)$.  Then $\displaystyle{\lim_{t\to +\infty}\|u(t)\|_X = 0}$.
\end{thm}
\noindent
This is proved for $X=\hdhalf(\rt)$ in \cite{GIP1}, and for $X=\ltrt$ in \cite{GIP3} (see also \cite{GIP2}).  The following result is classical (see, e.g., \cite{gallagher} for a proof):
\begin{thm}[Energy Inequality]\label{energysmall}
Let $u_0 \in \hdhalf(\rt)$ be any divergence-free field, and let $u=NS(u_0)$.  There exists a small universal constant $\d_0 >0$ such that $\|u_0\|_{\ltrt} \leq \d_0$ implies that $T^*(u_0) = +\infty$ and for any $t\geq 0$ the following energy inequality holds:
$$\|u(t)\|^2_{\hdhalf(\rt)} + \int_0^t \|\nabla u(s)\|^2_{\hdhalf(\rt)}\ ds \leq \|u_0\|^2_{\hdhalf(\rt)}\ .$$
\end{thm}
\noindent
Together with Theorem \ref{decay}, we now have the following:
\begin{cor}\label{energy}
Let $u_0\in \hdhalf(\rt)$, $\nabla \cdot u_0=0$ and $u=NS(u_0)$, and suppose $T^*(u_0)=+\infty$.  Then for some large $T_0>0$ we have, for all $t\geq T_0$,
$$\|u(t)\|^2_{\hdhalf(\rt)} + \int_{T_0}^t \|\nabla u(s)\|^2_{\hdhalf(\rt)}\ ds \leq \|u(T_0)\|^2_{\hdhalf(\rt)}\ .$$
\end{cor}
\noindent
Finally, we remark that the following fact follows easily from the local theory and Corollary \ref{energy} (recall $\se_T = L^\infty((0,T);\hdhalf) \cap L^2((0,T);\hdthalf)$\ ):
\begin{thm}\label{thm:blowup}
For $T\in (0,+\infty]$, take $(X_0,X_T)$ to be either $(\hdhalf(\rt),\se_T)$ or ${(\ltrt, L^5(\rt \times (0,T)))}$.  Then for any $u_0 \in X_0$,
$$T^*(u_0) < +\infty \quad \iff \quad \|NS(u_0)\|_{X_{T^*(u_0)}} = +\infty\ .$$
\end{thm}
\section{Main results}
The following ``blow-up" criterion for the Navier-Stokes equations was proved in \cite{ess}:
If $u$ is a {\em weak} ``Leray-Hopf" solution of the Navier-Stokes equations (essentially, a distributional solution to (\ref{eq:nseA}) satisfying a natural energy estimate) and if $(0,T^*)$ is the maximal interval on which $u$ is smooth, then if  $T^* < +\infty$, we must have
$$\limsup_{t\nearrow T^*} \|u(t)\|_3 = +\infty\ .$$
Alternatively, one can say that $\|u\|_{L^\infty((0,T^*);L^3)}<+\infty$ implies that $T^* = +\infty$.
\\\\
We would like to suggest a different approach to this statement, using the method of ``critical elements" introduced in \cite{km3,km1,km2}, in the setting of {\em mild} solutions.  We consider here a simpler special case of the above, namely the statement that for any $u_0\in \hdhalf$,
\begin{equation}\label{eq:s}
\|NS(u_0)\|_{L^\infty((0,T^*(u_0));\hdhalf)}<+\infty \quad \Longrightarrow \quad T^*(u_0) = +\infty\
\end{equation}
where $NS(u_0) \in \cap_{T<T^*(u_0)}\mathcal{E}^{(1/2)}_{T}$ (see (\ref{B})) is the associated mild solution.  Theorem \ref{decay} then implies moreover that  $\lim_{t\to +\infty}\|u(t)\|_{\hdhalf} =0$, so that (\ref{eq:s}) may be thought of as a statement of ``global existence and scattering" (time-asymptotic convergence to a solution of the linear equation), such as was considered in \cite{km3} for non-linear Schr\" odinger equations.  We start with a setup following \cite{km3} as follows:
\\\\
Suppose (\ref{eq:s}) is false; then there exists some maximal finite $A_c > 0$ (see (\ref{eq:an}) below) such that
$$ \|u\|_{L^\infty((0,T^*(u_0));\hdhalf)}<A_c \quad \Longrightarrow \quad T^* = +\infty\ $$
while, for any $A>A_c$, there exists some initial datum $u_{0,A}$ and solution $u_A$ with $T^*(u_{0,A}) < +\infty$ and $\|u_A\|_{L^\infty((0,T^*(u_{0,A}));\hdhalf)} \leq A$.  Note that the local theory implies $A_c >0$ is well-defined since $T^*(u_0) = +\infty$ (and moreover $u\in \mathcal{E}^{(1/2)}_{T}$ for $T=+\infty$ as well by Corollary \ref{energy}) for $\|u_0\|_{\hdhalf}$ small enough.
\\\\
Specifically, we define the critical value $A_c$ by
\begin{equation}\label{eq:an}
A_c := \sup \{A>0\ ;\  \|NS(u_0)\|_{L^\infty((0,T^*(u_0));\hdhalf)}\leq A  \Longrightarrow  T^*(u_0) = +\infty \}\ .
\end{equation}
Under the assumption that $A_c < +\infty$, we'll prove the following two theorems, whose analogs were used to prove statements similar to (\ref{eq:s}) in \cite{km3,km1,km2}:
\begin{thm}\label{thm:a}
Suppose $A_c < +\infty$.  Then there exists some $u_{0,c} \in \hdhalf$ with associated mild solution $u_c$ such that $T^*(u_{0,c}) < +\infty$ and $\|u_c\|_{L^\infty((0,T^*(u_{0,c}));\hdhalf)}=A_c$.
\end{thm}
\noindent
We call $u_c$ a ``critical element".
\begin{thm}\label{thm:b}
Suppose $A_c < +\infty$, and let $u_{0} \in \hdhalf$ with associated mild solution $u=NS(u_0)$ satisfy $T^*(u_{0}) < +\infty$ and $\|u\|_{L^\infty((0,T^*(u_{0}));\hdhalf)}=A_c$.  For any $\{t_n\} \subset [0,T^*(u_0))$ such that $t_n \nearrow T^*(u_0)$, there exist sequences $\{s_n\}$ with $t_n \leq s_n < T^*(u_0)$, $\{x_n\} \subset \rt$  and $\{\l_n\} \subset (0,\infty)$  with $\l_n \to +\infty$ such that a subsequence of $\frac{1}{\l_n} u\left( \frac{\ \cdot\  -\ x_n}{\l_n}, s_n    \right)$ converges in $L^3$.
\end{thm}
\noindent
This is a slightly weaker version of the compactness in $L^3$ of the closure of a family $K$, where
\begin{equation}\label{eq:famK}
K:= \left\{ \frac{1}{\l(t)} u\left( \frac{x-x(t)}{\l(t)}, t    \right) ,\ \ 0\leq t < T^*(u_0)         \right\}
\end{equation}
for some $\l(t)>0$, $x(t) \in \rt$, as was proved in \cite{km3,km1,km2}.  Note, however, that for any sequence $t_n \nrightarrow T^*(u_0)$, setting $x(t_n) \equiv 0$ and $\l(t_n) \equiv 1$, a convergent subsequence exists by the $L^3$-continuity of $u$.  The only remaining scenario is therefore treated by the weaker statement in Theorem \ref{thm:b}.
\\\\
The idea is that $A_c < +\infty$ implies the existence of a critical element which produces a type of compactness (see (\ref{eq:kay})) in the family $K$.  If one can rule out such compactness, then one has proved (\ref{eq:s}).  (The nonexistence of a particular example of a fully compact family K was established in \cite{necas} and \cite{tsai} regarding the so-called ``self-similar'' solutions.)  The program is therefore completed with the following ``ridigity" theorem:
\begin{thm}\label{thm:c}
Any $u_{0} \in \ltrt$ satisfying the conditions of Theorem \ref{thm:b} must be identically zero.
\end{thm}
\noindent
Theorem \ref{thm:c} has the immediate corollary that no such element $u_0$ can exist, since $T^*(0)=+\infty$ whereas by assumption $T^*(u_0)<+\infty$. This implies, due to Theorems \ref{thm:a} and \ref{thm:b}, that $A_c = +\infty$, which proves the regularity criterion (\ref{eq:s}) (Theorem \ref{thm:regularity}) as desired.
\subsection{Preliminary lemmas}
\noindent
In order to prove Theorem \ref{thm:a} (and also Theorem \ref{thm:b}), we will need the following ``profile decomposition" which was proved in \cite{gallagher}.  Recall the notation $NS(u_0)$ for the mild solution in $\cap_{T<T^*(u_0)} \mathcal{E}^{(1/2)}_{T}$ associated to an initial datum $u_0\in \hdhalf$.
\begin{thm}[Profile Decomposition]\label{thm:prof}
Let $\{u_{0,n}\}$ be a bounded sequence of divergence-free vectors in $\hdhalf$.  Then there exist $\{x_{j,n}\} \subset \rt$ and $\{\l_{j,n}\} \subset (0,+\infty)$ which are ``orthogonal" in the sense that, for $j,j' \in \N$, $j\neq j'$,
$$\textrm{either} \quad \lim_{n \to +\infty} \frac{\ljn}{\ljpn} + \frac{\ljpn}{\ljn}  = +\infty$$
$$\textrm{or} \quad \frac{\ljn}{\ljpn} \equiv 1  \quad \textrm{and} \quad \lim_{n \to +\infty}  \frac{|x_{j,n} - x_{j',n}|}{\ljn} = +\infty\ ,$$
and a sequence of divergence-free\footnote{Note that this property of each $V_{0,j}$ is a consequence of (\ref{eq:aa}) and the orthogonality of the scales-cores.  Indeed, $\tilde{u}_{0,n} \rightharpoonup V_{0,j}$ where $\tilde{u}_{0,n}(x) = \ljn u_{0,n}(\ljn x + \xjn)$ is divergence-free.} vector fields $\{V_{0,j}\} \subset \hdhalf$ with $T^*(\voj) < +\infty$ for an at most finite number of $j\in \mathbb{N}$ such that the following is true:  after possibly taking a subsequence in $n$,
\begin{equation}\label{eq:aa}
u_{0,n}(x) = \sum_{j=1}^{J} \frac{1}{\ljn} \voj \left( \frac{x-\xjn}{\ljn}\right) + w_n^J(x)\
\end{equation}
and
\begin{equation}\label{eq:ab}
u_{n}(x,t) = \sum_{j=1}^{J} \frac{1}{\ljn} U_j \left( \frac{x-\xjn}{\ljn}, \frac{t}{\ljn^2} \right) + w_n^{l,J}(x,t) + r_n^J(x,t)\
\end{equation}
for any $J\in \mathbb{N}$, for $x\in \rt$ and $t\in(0,T_n)$ where $T_n:= \min \{\ljn^2 T^j\ | \ T^*(\voj) < +\infty\}$ (or $T_n \equiv +\infty$ if $T^*(V_{0,j}) = +\infty$ for all $j$) for any fixed numbers $T^j < T^*(\voj)$, $u_n = NS(u_{0,n})$, $U_j = NS(V_{0,j})$ and $w_n^{l,J}(t) = \etl w_n^J$ for some $w_n^J \in \hdhalf$ and $r_n^J \in \E_{T_n}$, and the following orthogonality properties hold:
\begin{equation}\label{eq:ac}
\|u_{0,n}\|_{\hdhalf}^2 = \sum_{j=1}^{J} \| \voj \|_{\hdhalf}^2+ \|w_n^J\|_{\hdhalf}^2 + \e^J_n \ , \quad \lim_{n\to +\infty}\e^J_n = 0 \quad \forall  J\in \N,
\end{equation}
\begin{equation}\label{eq:ad}
\lim_{J\to +\infty} \left( \limsup_{n\to +\infty} \|w_n^J\|_3 \right) = 0\ ,
\end{equation}
and
\begin{equation}\label{eq:ae}
\lim_{J\to +\infty} \left( \limsup_{n\to +\infty} \|r_n^J\|_{\E_{T_n}} \right) = 0\ .
\end{equation}
\end{thm}
\noindent
Recall that $\et$ is defined for $T>0$ by (\ref{eq:defet}), and note that $T_n$ is simply some time such that for $t\in (0,T_n)$, (\ref{eq:ab}) avoids the finite blow-up times of any of the $\voj$'s due to the natural scaling of the equation: if $v$ is a solution on $(0,T)$ then $v_\l(x,t) := \l v(\l x, \l^2 t)$ is a solution on $(0,\frac{T}{\l^2})$ for any $\l >0$.
\\\\
We'll also need the following fact regarding the ``profiles" $U_j$ which appear in Theorem \ref{thm:prof}:
\begin{clm}\label{lemma:aa}
If $T^*(\voj) = +\infty$ for all $j \geq 1$, then there exists some $n_0 \geq 1$ such that $T^*(u_{0,n_0}) = +\infty$.
\end{clm}
\noindent
This is in fact implicit in the statement of Theorem \ref{thm:prof}.  The proof can be found in \cite{gallagher}, for example if one follows the
proof of Theorem 2 part (iii) in section 3.2.2 of that paper.  One sees that,
under the assumptions of Claim 3.5 above, the remainder $r_n^J$ belongs to
$E_{T^*(u_{0,n})}$ for large $J$ and $n$ and hence $u_n$ belongs to the same
space due to (\ref{eq:ab}) and standard heat estimates.  Theorem 2.12 then implies that $T^*(u_{0,n}) = +\infty$ for such $n$. \hfill $\Box$
\\\\
Finally, we will need the following ``backwards uniqueness" type of lemma (a stronger version of which is proved in the last section):
\begin{lemma}\label{lemma:zero}
Let $u(t)$ be a mild solution of NSE with initial datum $u(0) = u_0 \in \lt$.  Suppose $u(t_1) = 0$ for some $t_1 \in (0,T^*(u_0))$.  Then $u\equiv 0$. In particular, $u_0 = 0$ and $T^*(u_0) = +\infty$.
\end{lemma}
\noindent
Lemma \ref{lemma:zero} will follow from the following backwards uniqueness theorem for systems of parabolic equations, which was proved in \cite{ess2,ess3} (for the more general situation of a half-space see also \cite{ess}, but we do not require such generality):
\begin{thm}[Backwards Uniqueness]\label{thm:backwards}
Fix any $R, \d, M$ and $c_0 >0$.  Let $Q_{R,\d}:= (\rt \backslash B_R(0)) \times (-\delta,0)$, and suppose a vector-valued function $v$ and its distributional derivatives satisfy $v,\ v_t,\ \nabla v, \ \nabla^2 v \in L^2(\Omega)$ for any bounded subset $\Omega \subset Q_{R,\d}$, $|v(x,t)| \leq e^{M|x|^2}$ for all $(x,t)\in Q_{R,\d}$, $|v_t - \Delta v| \leq c_0(|\nabla v| + |v|)$ on $Q_{R,\d}$ and $v(x,0) = 0$ for all $x\in\rt \backslash B_R(0)$.  Then $v\equiv 0$ in $Q_{R,\d}$.
\end{thm}
\noindent
{\bf Proof of Lemma \ref{lemma:zero}:}
\quad
We will apply Theorem \ref{thm:backwards} with $v=\omega = \textrm{curl}_xu$, the associated vorticity to the velocity field $u$.  Suppose for the moment that $\omega$ satisfies the assumptions of Theorem \ref{thm:backwards} (in the whole space, in fact).  Then $u(t_1) = 0$ implies $\omega(t_1) = 0$, and hence $\omega \equiv 0$ in $\rt \times (t_1 - \delta,t_1)$, for any $\d \in (0,t_1)$.  Since mild solutions to NSE satisfy $\textrm{div}_xu = 0$ in the distributional sense, we see that $\Delta_xu(t) = 0$, $u(t) \in L^3$ for each $t\in (t_1 - \delta,t_1)$ which implies that $u \equiv 0$ in $\rt \times (t_1 - \d,t_1)$.  Taking $\d \approx t_1$ and using continuity at $t=0$, the theorem follows. (In the forward direction, one uses the uniqueness of Theorem \ref{unique}.)
\\\\
The vorticity satisfies the equation
$$\omega_t - \Delta \omega + (u\cdot \nabla) \omega - (\omega \cdot \nabla) u = 0\ $$
in the sense of distributions.  The assumptions of Theorem \ref{thm:backwards} follow therefore from the fact that $u\in L^\infty (\rt \times (\e,T^*(u_0) - \e))$ for any $\e > 0$, which follows from Remark \ref{remark:bdd} and arguments similar to those in the proof of Lemma \ref{lemma:a}.  The bounds on the derivatives of $u$ (and hence of $\omega$) are then standard, see for example \cite{knss}. \hfill $\Box$
\subsection{Existence of a critical element}
\noindent
In this section we'll prove Theorem \ref{thm:a}, which establishes the existence of a critical element.
\\\\
{\bf Proof of Theorem \ref{thm:a}:}
\quad
Define $A_c$ by (\ref{eq:an})
and assume that $A_c < +\infty$.  (Note that $A_c$ is well-defined by global existence for small data.)
By definition of $A_c$, we can pick $A_n \searrow A_c$ as $n \to +\infty$ and $u_{0,n}\in \hdhalf$, $\nabla \cdot u_{0,n} = 0$ such that
\begin{equation}\label{eq:ah}
T^*(u_{0,n}) < +\infty \quad \textrm{and} \quad \|NS(u_{0,n})\|_{L^\infty((0,T^*(u_{0,n}));\hdhalf)}\leq A_n
\end{equation}
holds for all $n$.  One cannot prove, as one would naturally hope, that a subsequence of these $u_{0,n}$ converges to some $u_{0,c}$ satisfying the assertion of the theorem, but we will see that something similar is true. We may assume that
\begin{equation}\label{eq:af}
\|\uon\|_{\hdhalf} \leq A_n \leq 2A_c\ ,
\end{equation}
and hence we may apply the profile decomposition (Theorem \ref{thm:prof}) to the sequence $\{\uon\}$.  For the remainder of this section, $\voj,\ U_j,\ u_n$ etc. will refer to this particular sequence.  We will complete the proof of the theorem by showing that $V_{0,j_0}$ satisfies the assertion of the theorem for some $j_0$, and the profile $U_{j_0}=NS(V_{0,j_0})$ will be our critical element.
\\\\
Define $T^*_j:= T^*(\voj)$, and note that (\ref{eq:af}) and (\ref{eq:ac}) imply that
\begin{equation}\label{eq:ag}
\sum_{j=1}^{\infty} \|\voj\|^2_\hdhalf < +\infty\ ,
\end{equation}
so that $\lim_{j\to \infty} \|\voj\|_\hdhalf = 0$ and therefore by the small data result, $T^*_j = +\infty$ for large enough $j$.  Corollary \ref{energy} then implies that $\|U_j\|_{\E_{(+\infty)}}  <+\infty$ for large $j$ as well.
\\\\
Note that Theorem \ref{thm:blowup} and Claim \ref{lemma:aa} together imply that there exists some $j\geq 1$ such that $\|U_j\|_{\E_{T^*_j}} = +\infty$, and hence we may re-order the profiles in the decomposition such that for some $J_1 \geq 1$,
\begin{equation}\label{eq:ak}
\|U_j\|_\etj = +\infty  \ \ \textrm{for} \ \  1\leq j \leq J_1\
\end{equation}
and
$$  \|U_j\|_\etj< +\infty \ \  \textrm{for} \ \  j > J_1 \ .$$
Note also that Theorem \ref{thm:blowup} and Corollary \ref{energy} imply that
\begin{equation}\label{eq:al}
T^*_j < +\infty \quad  \textrm{for} \quad 1\leq j \leq J_1 \qquad \textrm{and} \qquad  T^*_j = +\infty \quad  \textrm{for}\quad j > J_1\ .
\end{equation}
To prove Theorem \ref{thm:a}, it suffices therefore to show that
\begin{equation}\label{eq:ai}
\|U_{j_0}\|_{L^\infty((0,T^*_{j_0});\hdhalf)} = A_c
\end{equation}
for some $j_0\in \{1, \dots , J_1\}$.
This will be accomplished using the following claim, which extends the orthogonality property (\ref{eq:ac}):
\begin{clm}\label{clm:d}
Set $T^*_{j,k} := T^*_j - \frac{1}{k}$ and $t^n_k := \min_{1\leq j \leq J_1} \ljn^2 T^*_{j,k}$ for $j, k \in \N$.  Fix some $k\geq 1$, and let $\{t_n\} \subset \R$ be any sequence such that $t_n \in [0,t^n_k]$ for all $n$. Then there exist subsequences $n(m,k)$ and $J(m,k)$ depending on $k$ and indexed by $m$ such that for $n=n(m,k)$ and $J=J(m,k)$,
$$\|u_n(t_n)\|_\hdhalf^2 = \sum_{j=1}^J \|\tilde{U}_{j,n}(t_n)\|_\hdhalf^2 + \|w_n^{l,J}(t_n)\|_\hdhalf^2 + \e(k,m)\ ,$$
where $\displaystyle{\lim_{m\to +\infty}\e(k,m) = 0}$ for each $k$ and we have set $\tilde{U}_{j,n}(x,t) := \frac{1}{\ljn}U_j \left( \frac{x- \xjn}{\ljn}, \frac{t}{\ljn^2} \right)$.
\end{clm}
\noindent
Assuming Claim \ref{clm:d} momentarily, we prove (\ref{eq:ai}) (and hence Theorem \ref{thm:a}) as follows:
\\\\
Let $t^n_{j,k} := \ljn^2 T^*_{j,k}$, so that $t^n_k = \min_{1\leq j \leq J_1} t^n_{j,k}$.  For any fixed $k$, there exists some $j_0^k \in \{1, \dots , J_1\}$ such that
$$t^n_k = t^n_{j^k_0,k}$$
for infinitely many $n$.  Also, there exists some $j_0\in \{1, \dots , J_1\}$ such that $j^k_0 = j_0$ for infinitely many $k$, say for a subsequence $k_\alpha\to \infty$ as $\alpha \to +\infty$.
For each fixed $\alpha$, take $m_{\alpha,\beta}$ to be a subsequence of $m$'s indexed by $\beta \to +\infty$ such that
\begin{equation}\label{eq:aj}
t^{n(m_{\alpha ,\beta},k_\alpha)}_{k_\alpha} = t^{n(m_{\alpha ,\beta},k_\alpha)}_{j_0,k_\alpha}\
\end{equation}
(where $n(m,k)$, $J(m,k)$ are as in Claim \ref{clm:d}), and $m_{\alpha,\beta} \to +\infty$ as $\beta \to +\infty$.
\\\\
Since $1\leq j_0 \leq J_1$, we have $\|U_{j_0}\|_{\E_{T^*_{j_0}}} = +\infty$ by (\ref{eq:ak}).
Let
$$A^{(j_0)}:= \sup_{0\leq t \leq T^*_{j_0}} \|U_{j_0}(t)\|_\hdhalf
\quad \textrm{and} \quad
A^{(j_0)}_k := \sup_{0\leq t \leq T^*_{j_0,k}} \|U_{j_0}(t)\|_\hdhalf$$
so that $A^{(j_0)}_k \nearrow A^{(j_0)}$ as $k\to +\infty$.  Note that $A^{(j_0)}\geq A_c$, else we would have $T^*_{j_0} = +\infty$ by the definition of $A_c$ which would contradict (\ref{eq:al}).
Theorem \ref{thm:a} will be proved if we show $A^{(j_0)} = A_c$ (whereby we can set $u_{0,c} = V_{0,j_0},\ u_c = U_{j_0}$), so it remains only to show that $A^{(j_0)} \leq A_c$, which we will now do.
\\\\
By continuity in $\hdhalf$, we may take $t_k \in [0,T^*_{j_0,k}]$ such that $\|U_{j_0}(t_k)\|_\hdhalf = A^{(j_0)}_k$.  Set
$t_{k,n} = \l^2_{j_0,n}t_k$, and note that for $k=k_\alpha$ and $n=n(m_{\alpha,\beta},k_\alpha)$ we have
$$0 \leq t_{k,n} \leq \l^2_{j_0,n} T^*_{j_0,k} = t_k^n$$
by (\ref{eq:aj}).  We may therefore apply Claim \ref{clm:d} with $t_n = t_{k,n} = t_{k_\alpha,n(m_{\alpha,\beta},k_\alpha)}$ to conclude that for fixed $\alpha$, there exists a subsequence of $\beta$'s such that for $n=n(m_{\alpha,\beta},k_\alpha)$, $J=J(m_{\alpha,\beta},k_\alpha)$ and $k=k_\alpha$ we have
\begin{equation}\label{eq:am}
\|u_n(t_{k,n})\|_\hdhalf^2 = \sum_{j=1}^J \|\tilde{U}_{j,n}(t_{k,n})\|_\hdhalf^2 + \|w_n^{l,J}(t_{k,n})\|_\hdhalf^2 + \e(\alpha,\beta)\ ,
\end{equation}
where
$\e(\alpha,\beta) \to 0$ as $\beta \to +\infty$.  Recall from (\ref{eq:ah}) that we defined $A_n$ so that
$ \sup_{0\leq t \leq T^*(u_{0,n})} \|u_n(t)\|_\hdhalf \leq A_n.$
Therefore, by (\ref{eq:am}), we have
$$A^2_{n(m_{\alpha,\beta},k_\alpha)} \geq (A^{(j_0)}_{k_\alpha})^2 + \e(\alpha,\beta)\ .$$
Keeping $\alpha$ fixed and letting $\beta \to +\infty$, this gives $A_c^2 \geq (A^{(j_0)}_{k_\alpha})^2$.  Finally, letting $\alpha \to +\infty$, we see $A_c^2 \geq (A^{(j_0)})^2$ as desired proving Theorem \ref{thm:a}. \hfill $\Box$
\begin{remark}\label{remark:a}
These arguments can moreover be used to show that in fact $j_0 = J_1 = 1$.  That is, there is one and only one profile $U_1$ which blows up in finite time.  This is due to the fact that we have now shown that $A^{(j_0)}_{k_\alpha} \to A_c$ as $\alpha \to +\infty$, hence (\ref{eq:am}) can be used to show that $\|\tilde{U}_{j,n}(t_{k,n})\|_\hdhalf$ can be made arbitrarily small for any $j\neq j_0$, which implies that $T^*_j = +\infty$ by the global existence for small data result and the natural scaling of the equations.
\end{remark}
\ \\
\textbf{Proof of Claim \ref{clm:d}:}
\quad
For $J\in \N$, define
$$\tU_{J,n}:= \soj \tU_{j,n}\qquad \textrm{and} \qquad \tU^J_{J_2,n}:= \tU_{J,n} - \tU_{J_2,n} \quad \textrm{for} \quad 1\leq J_2 \leq J\ .$$
Fix some $k\in \N$.  We now claim that the following three statements hold, for some subsequences $n=n(m)$ and $J=J(m)$ indexed by m:
\begin{enumerate}
\item For any $j,j' \in \N$, $j\neq j'$,
\begin{equation}\label{eq:ddd}
< \doh \tU_{j,n}(t_n),\doh \tU_{j',n}(t_n) > \ \ \longrightarrow\ \  0 \quad \textrm{as} \quad m\to \infty\ .
\end{equation}
\item For any $J \geq J_2 \geq 1$,
\begin{equation}\label{eq:dde}
\|\tU^J_{J_2,n}(t_n)\|^2_\hdhalf \leq 2 \sum_{j=J_2+1}^J \|\tU_{j,n}(t_n)\|^2_\hdhalf \quad \textrm{for all} \ m \ .
\end{equation}
\item There exists some $M>0$ such that
\begin{equation}\label{eq:ddf}
\|\tU_{J,n}(t_n)\|_\hdhalf \leq M \quad \textrm{for all} \ m \ .
\end{equation}
\end{enumerate}
The easiest to show is (\ref{eq:ddf}):  Since $\tujn(t_n)$ is defined for all $j\geq 1$ by assumption, the properties of the profile decomposition imply\footnote{This fact is used again below in Claim \ref{clm:e}, see (\ref{eq:dda}).} (see \cite{gallagher}) that $t_n < T^*(\uon)$ for all $n$, and so $\|u_n(t_n)\|_\hdhalf \leq 2A_c$ by (\ref{eq:ah}).  Using (\ref{eq:ac}) and heat estimates, we see
$$\|w_n^{l,J}(t_n)\|_\hdhalf \leq \|w_n^J\|_\hdhalf \leq 3A_c$$
for large enough $n$.  For each $m$, use (\ref{eq:ae}) to take $J$ large enough, and a corresponding sufficiently large $n$ (both in an increasing fashion) such that
$$\sup_{t\in[0,t_k^n]}\|r^J_n(t)\|_\hdhalf \leq 1\ ,$$
whereupon (\ref{eq:ab}) gives
$\|\tU_{J,n}(t_n)\|_\hdhalf \leq 2A_c + 3A_c + 1 =:M$.
\\\\
To prove (\ref{eq:ddd}), we take the following diagonal-type subsequence:  If $t_{1,n}:={t_n/\l_{1,n}}$ is bounded, pass to a subsequence such that $t_{1,n} \to t_1 \in [0,+\infty)$.  Take $n(1,k)$ to be the first element of this subsequence (or of the original subsequence if $t_{1,n}$ is unbounded).  (This is $m=1$.)
Assume now that the first $m-1$ values of $n$ are chosen.
If $t_{m,n}:={t_n/\l_{m,n}}$ is bounded, pass to a subsequence of the $(m-1)$st subsequence such that $t_{m,n} \to t_m \in [0,+\infty)$.  Take $n(m,k)$ to be the first element of this sequence which comes after $n(m-1,k)$ (or of the original $(m-1)$st subsequence if $t_{m,n}$ is unbounded).  In this way, $t_{j,n} \to t_j \in [0,+\infty)$ as $m\to \infty$ ($n=n(m,k)$) for any $j$ such that $t_{j,n}$ is bounded.
\\\\
Note that for $1\leq j\leq J_1$, $t_{j,n} = {t_n/\l_{j,n}} \leq T^*_{j,k} < +\infty$.  Therefore if $t_{j,n} \to +\infty$, necessarily $j>J_1$ and $U_j \in \se_\infty$.  Hence by Theorem \ref{decay}, $\|U_j(t)\|_\hdhalf \to 0$ as $t\to +\infty$, and hence
\begin{equation}\label{eq:ddg}
\|\tU_{j,n}(t_{j,n})\|_\hdhalf = \|U_j (t_{n})\|_\hdhalf \ \ \longrightarrow \ \ 0 \quad \textrm{as} \quad n\to \infty\ .
\end{equation}
In the expression in (\ref{eq:ddd}), whenever $t_{j,n}$ is bounded, using the $\hdhalf$-continuity of $U_j$ in a neighborhood of $t_j$, we can replace
$$\doh \left\{\frac{1}{\ljn}U_j\left( \frac{x-\xjn}{\ljn},t_{j,n}\right)\right\} \qquad \textrm{by} \qquad
\doh \left\{\frac{1}{\ljn}U_j\left( \frac{x-\xjn}{\ljn},t_j\right)\right\}\ .$$
We can also assume that $\doh U_j(t_j) \in \mathcal{C}^\infty_0$ by approximating in $L^2_x$.
If both $t_{j,n}$ and $t_{j',n}$ are unbounded, then
$$|<\doh \tU_{j,n}(t_n),\doh \tU_{j',n}(t_n)>| \leq \|\tU_{j,n}(t_n)\|_\hdhalf  \|\tU_{j',n}(t_n)\|_\hdhalf =$$$$= \|U_j(t_{j,n})\|_\hdhalf \|U_j(t_{j',n})\|_\hdhalf \longrightarrow 0 \quad \textrm{as} \quad n\to \infty$$
by (\ref{eq:ddg}).
If $\tjn$ is bounded and $\tjpn$ is unbounded, we estimate the above expression by
$$\|U_j(t_j)\|_\hdhalf \|U_{j'}(\tjpn)\|_\hdhalf \to 0$$
as $n\to \infty$ since we have replaced $t_{j,n}$ by the constant $t_j$ in the first term, and the second term tends to zero.
If both $\tjn$ and $\tjpn$ are bounded, replace $\tjn$ by $t_j$ and $\tjpn$ by $t_{j'}$, and note by a change of variables that
$$\left<\doh \left\{\frac{1}{\ljn}U_j\left( \frac{x-\xjn}{\ljn},t_j\right)\right\},\doh \left\{\frac{1}{\l_{j',n}}U_{j'}\left( \frac{x-x_{j',n}}{\l_{j',n}},t_{j'}\right)\right\}\right> =$$
$$=\left(\frac{\ljn}{\ljpn}\right)^\frac{3}{2} \int \doh U_j(y,t_j)\doh U_{j'}\left(\frac{\ljn}{\ljpn}y + \frac{\xjn - \xjpn}{\ljpn},t_{j'}\right)\ dy$$
which, since we've assumed the functions in the integrand are compactly supported, is easily seen to tend to zero as $n\to \infty$ if ${\ljn/\ljpn} \to 0$ or if $\ljn = \ljpn$ and $|{(\xjn-\xjpn)/\ljpn}| \to +\infty$.  The case ${\ljpn/\ljn} \to 0$ is handled similarly, and (\ref{eq:ddd}) is proved.
\\\\
(\ref{eq:dde}) is proved using (\ref{eq:ddd}) as follows:  Note that
$$\|\tU^J_{J_2,n}(t_n)\|^2_\hdhalf = \sum_{j=J_2+1}^J \|\tujn(t_n)\|^2_\hdhalf + \sum_{J_2 \leq j\neq j' \leq J} <\doh \tujn(t_n),\doh \tujpn (t_n)>\ .$$
Suppose the sum on the right has $C(J_2,J)$ terms in it.  For a fixed $J$, (\ref{eq:ddd}) allows us to take $n$ large enough (that is, we take a further subsequence of $n(m,k)$\ ) so that all the terms are smaller than ${\e/C(J_2,J)}$, where $\e = \sum_{j=J_2+1}^J \|\tujn(t_n)\|^2_\hdhalf$.  As $J$ increases with $m$, we let $n$ increase sufficiently rapidly with $m$ as well so that this is true for all $m$, and (\ref{eq:dde}) is proved.
\\\\
Returning now to the proof of Claim \ref{clm:d}, we want to show that for any $\e>0$, there exists some $M_1 >0$ such that, for all $m\geq M_1$, $|A|<\e$ where
$$A:= \|u_n(t_n)\|^2_\hdhalf - \soj \|\tujn(t_n)\|^2_\hdhalf - \|w_n^{l,J}(t_n)\|_\hdhalf^2\ .$$
We write $u_n = \tU_{J_2,n} + \tU^J_{J_2,n} + \wnlj + \rnj$ for a sufficiently large $J_2$ to be chosen, expand the first term in $A$ and cancel the equal terms.  We'll then split $A$ into $I + II$, where $II$ contains all terms with $\tujn$ for $j>J_2$ (see below).  Note that for $j> J_2 > J_1$, $U_j \in \se_\infty$ and hence, for $J_2$ sufficiently large, Theorem \ref{energysmall} and (\ref{eq:ag}) give the energy estimate $\|U_j\|_{\se_\infty} \lesssim \|\voj\|_\hdhalf$.  Then, using (\ref{eq:ac}), (\ref{eq:ag}) and (\ref{eq:dde}), for any $\tilde{\e}>0$ we have
\begin{equation}\label{eq:ddj}
\|\tU^J_{J_2,n}(t_n)\|^2_\hdhalf \leq 2 \sum_{j=J_2+1}^J \|\tujn(t_n)\|^2_\hdhalf \lesssim \sum_{j=J_2+1}^J \|\voj\|^2_\hdhalf < \tilde{\e}
\end{equation}
for $J_2 = J_2(\tilde{\e})$ sufficiently large.  Now let
$$II:= 2<\doh \tU^J_{J_2,n}(t_n), \doh \{\tU_{J,n}(t_n) + \wnlj (t_n) + \rnj (t_n)\} > +$$$$+ \|\tU^J_{J_2,n}(t_n)\|^2_\hdhalf - \sum_{j=J_2 +1}^J \|\tujn (t)\|^2_\hdhalf\ .$$
We estimate
$$|II| \leq 2\|\tU^J_{J_2,n}(t_n)\|_\hdhalf \cdot \left( \|\tU_{J,n}(t_n)\|_\hdhalf + \|\wnlj(t_n)\|_\hdhalf + \|\rnj(t_n)\|_\hdhalf \right) + $$$$+ \|\tU^J_{J_2,n}(t_n)\|^2_\hdhalf + \sum_{j=J_2 +1}^J \|\tujn (t)\|^2_\hdhalf\ .
$$
Noting that we have
$$\|\tU_{J,n}(t_n)\|_\hdhalf  \leq M\ , \quad  \|\wnlj(t_n)\|_\hdhalf \leq \|\wnj\|_\hdhalf \leq 3A_c \quad \textrm{and} \quad \|\rnj(t_n)\|_\hdhalf \leq 1$$
by (\ref{eq:ddf}), heat estimates and (\ref{eq:ac}), and by our choice of subsequences, respectively.  Using these and (\ref{eq:ddj}), for any $\e >0$ we can make $|II|<{\e/2}$ for large enough $J_2$, which we now fix.  We are now left to estimate
$$I:= \|\tU_{J_2,n}(t_n) + \wnlj (t_n) + \rnj(t_n)\|^2_\hdhalf - \sum_{j=1}^{J_2} \|\tujn (t_n)\|^2_\hdhalf - \|\wnlj(t_n)\|^2_\hdhalf$$
$$= \sum_{1\leq j\neq j' \leq J_2} <\doh \tujn (t_n),\doh \tujpn(t_n)> +\qquad \qquad$$$$ \qquad \qquad + 2<\doh \tU_{J_2,n}(t_nn), \doh \{ \wnlj (t_n) + \rnj (t_n) \}> $$$$\qquad \qquad + 2<\wnlj (t_n),\rnj(t_n)> + \|\rnj(t_n)\|^2_\hdhalf\ .$$
The first term is handled by (\ref{eq:ddd}).  Since we may pass to a subsequence such that $\|\rnj(t_n)\|_\hdhalf \leq \|\rnj\|_{\se_{t_k^n}} \to 0$ as $m \to \infty$ by (\ref{eq:ae}), using (\ref{eq:ddf}) and the fact that $\|\wnlj(t_n)\|\leq 3A_c$ we see that the terms with $\rnj$ are small for large $m$.  The only remaining issue is
$$<\doh \tU_{J_2,n}(t_n),\doh \wnlj(t_n)>\ .$$
Consider $<\doh \tujn(t_n),\doh \wnlj(t_n)>$ for $1\leq j \leq J_2$.  If $\tjn \to +\infty$, then, since $\|\wnlj(t_n)\|_\hdhalf \leq 3A_c$, we can use again the fact that $\|U_j(\tjn)\|_\hdhalf \to 0$ as $n\to \infty$ to see that the term is small for large $m$ (recall that $n=n(m)$).  Otherwise, it suffices to consider (by replacing $\tjn$ by $t_j$ as before)
$$\left<\doh \left\{\frac{1}{\ljn}U_j\left(\frac{x-\xjn}{\ljn},t_j\right)\right\},\doh \wnlj(t_n)\right>=$$$$=
\int \doh U_j(y,t_j)\doh e^{\tjn \D}[\ljn \wnj (\ljn \cdot + \xjn)](y)\ dy\ .$$
Define $h_n(y):= e^{\tjn \D}[\ljn \wnj (\ljn \cdot + \xjn)](y)$.  We claim that $h_n \rightharpoonup 0$ in $\hdhalf$ and hence this term is small for large $m$ as well (for each $j\in [1, \dots , J_2]$), and Claim \ref{clm:d} is proved.  Indeed, note that $\|h_n\|_\hdhalf \leq 3A_c$ hence $h_n$ has a weak accumulation point in $\hdhalf$.  Note also that $\|h_n\|_3 \lesssim \|\wnj\|_3 \to 0$ as $m\to \infty$ (after passing to a subsequence) by (\ref{eq:ad}), so in particular $h_n \rightharpoonup 0$ in $\lt$.  Any weak limit in $\hdhalf$ for any subsequence is therefore also zero due to the embedding $(L^3)' \hookrightarrow (\hdhalf)'$, hence zero is the only weak accumulation point of $h_n$ in $\hdhalf$ and thus $h_n \rightharpoonup 0$ in $\hdhalf$.   This concludes the proof of Claim \ref{clm:d} (and Theorem \ref{thm:a}). \hfill $\Box$
\subsection{Compactness of critical elements}
\noindent
In this section, we'll prove Theorem \ref{thm:b} which establishes the compactness of any critical element arising from Theorem \ref{thm:a} under the assumption that $A_c < +\infty$.
\\\\
{\bf Proof of Theorem \ref{thm:b}:}
\quad
We will need the following claim:
\begin{clm}\label{clm:e}
Define $A_c$ by (\ref{eq:an}), and suppose that $A_c < +\infty$.  Suppose $u_0\in\hdhalf$ and $u=NS(u_0)$ are such that $T^*(u_{0}) < +\infty$ and $\|u\|_{L^\infty((0,T^*(u_{0}));\hdhalf)}=A_c$.   For any sequence $\{t_n\}$ such that $t_n \nearrow T^*(u_{0})$, let $\voj$, $U_j$ be the profiles associated with the sequence $u_{0,n} := u(t_n)$.  Then, after re-ordering, $T^*(\voj) < +\infty \iff j=1$, and
\begin{equation}\label{eq:dda}
T^*(u_0) - t_n \geq \l_{1,n}^2 T^*(V_{0,1})
\end{equation}
 for all $n$.  Moreover, after passing to a subsequence in $n$, for $j\geq 2$ either $U_j \equiv 0$ or
\begin{equation}\label{eq:as}
\frac{\l_{1,n}}{\ljn} \to +\infty \quad \textrm{as} \quad n\to +\infty\ .
\end{equation}
\end{clm}
\noindent
\textbf{Proof of Claim \ref{clm:e}:}
\quad
First note that the following considerations hold for {\em any} sequence $\{t_n\} \subset [0,T^*(u_0))$:
\\\\
Letting $\uon = u(t_n)$, $u_n(t) = NS(\uon)(t) = u(t_n + t)$ and noting $T^*(u_0) < +\infty$ implies that $\|u_n\|_{\se_{T^*(\uon)}} = +\infty$ for all $n\in \N$ by Theorem \ref{thm:blowup}, we can think of $\uon$ as a minimizing sequence for $A_c$, since $\|u(t_n)\|_\hdhalf \leq A_c \equiv: A_n$ for all $n$.  Proceeding as before and applying the profile decomposition, we see that, for a subsequence,
$$u(t_n) = \soj \frac{1}{\ljn}\voj \left( \frac{\cdot-\xjn}{\ljn} \right) + w_n^J(\cdot)$$
and, by uniqueness,
$$u(t_n+t) = u_n(t) =  \soj \frac{1}{\ljn}U_j \left( \frac{\cdot-\xjn}{\ljn}, \frac{t}{\ljn^2} \right) + w_n^{l,J}(\cdot, t) + r_n^J(\cdot,t)$$
for $t\in [0,t^n_k]$, where $t_k^n:= \l_{1,n}^2 \left(T^*(U_1) - \frac{1}{k}\right)$ for $k\in \N$ with the orthogonality properties (\ref{eq:ac}) - (\ref{eq:ae}).  This is justified since we saw before (see (\ref{eq:al}), (\ref{eq:ai}), Remark \ref{remark:a}, etc.) that we may re-order the elements so that $T^*(V_{0,1}) < +\infty$, $T^*(V_{0,j}) = +\infty$ for all $j\geq 2$ and
$$\|U_1\|_{\se_{T^*(V_{0,1})}} = +\infty\ , \qquad \|U_j\|_{\se_\infty} < +\infty \quad \textrm{for all} \quad j\geq 2\ ,$$
$$\sup_{t\in [0,T^*(\voo))}\|U_1(t)\|_\hdhalf = A_c$$
and it is noted moreover in \cite{gallagher} that equation (\ref{eq:dda}) holds for all $n$.  (In fact, (\ref{eq:dda}) is a simple consequence of the properties of the profile decomposition.)
\\\\
Suppose now that $t_n \nearrow T^*(u_0)$.  Then (\ref{eq:dda}) gives
\begin{equation}\label{eq:aaaaa}
\l^2_{1,n} \leq \frac{T^*(u_0) - t_n}{T^*(\voo)} \longrightarrow 0 \quad \textrm{as} \quad n\to \infty\ .
\end{equation}
We will now show that moreover, for $j\geq 2$, either $U_j \equiv 0$ or, after possibly passing to a subsequence, the limit (\ref{eq:as}) holds
as follows:  Take $t_k \in [0,T^*(\voo) - \frac{1}{k}]$ such that
$$\|U_1(t_k)\|_\hdhalf = \sup_{0\leq t \leq T^*(\voo)-\frac{1}{k}} \|U_1(t)\|_\hdhalf =: A_k \nearrow A_c\ .$$
Then, letting $t_{k,n}:= \l^2_{1,n}t_k$, we can apply Claim \ref{clm:d} (where actually $J_1 = 1$) to see that there exist subsequences $n=n(m,k)$, $J=J(m,k)$ indexed by $m$ for fixed $k$ such that
$$\|u(t_n + t_{k,n})\|_\hdhalf^2 = \soj \|\tujn(t_{k,n})\|_\hdhalf^2 + \|w_n^{l,J}(t_{k,n})\|_\hdhalf^2 + \e(k,m)$$
where $\e(k,m) \to 0$ as $m\to \infty$ for fixed $k$.  Therefore, as in Remark \ref{remark:a}, one can take a subsequence $m=m(k)$ such that, for $j\geq 2$,
$$\left\|U_j\left(\frac{\l^2_{1,n(m(k),k)}}{\l^2_{j,n(m(k),k)}}t_k\right) \right\|_\hdhalf
=\|\tilde{U}_{j,n(m(k),k)}(t_{k,n(m(k),k)})\|_\hdhalf \longrightarrow 0 \quad \textrm{as} \quad k\to \infty\ .
$$
If along this subsequence $n=n(m(k),k)$ we had that ${\l^2_{1,n}t_k/\ljn^2}$ were bounded, one could therefore extract a convergent subsequence with limit $\bar{t} < \infty$, and therefore conclude that $U_j(\bar{t}) = 0$ by the $\hdhalf$-continuity and hence that $U_j \equiv 0$ by Lemma \ref{lemma:zero}.  Since $0\leq t_k < T^*(\voo) < \infty$, if $U_j \neq 0$ then the limit (\ref{eq:as}) must therefore hold for some subsequence, and Claim \ref{clm:e} is proved.
 \hfill $\Box$
\\\\
Returning now to the proof of Theorem \ref{thm:b}, fix a sequence $\{t_n\}$ such that $t_n \nearrow T^*$.  We define\footnote{ The specific factor ${1/2}$ in the definition of $s_n$ is chosen only for simplicity;  one may in fact replace ${1/2}$ by $\eta$ for any fixed $\eta \in (0,1)$.} a corresponding sequence $\{s_n\}$ by
\begin{equation}\label{eq:aq}
\frac{s_n - t_n}{\l^2_{1,n}} = \tfrac{1}{2}T^*(V_{0,1})
\end{equation}
for each $n$.  Note that, for $j\geq 2$, (\ref{eq:as}) gives (after passing to a subsequence in $n$)
\begin{equation}\label{eq:ar}
\frac{s_n - t_n}{\l^2_{j,n}} = \frac{T^*(V_{0,1})}{2}\frac{\l^2_{1,n}}{\ljn^2} \to + \infty \quad \textrm{as} \quad n\to +\infty\ .
\end{equation}
Note also that (\ref{eq:dda}) implies that $s_n \in (t_n,T^*(u_0))$.  Indeed, we have
$$0< t_n < s_n = t_n + \tfrac{1}{2}\l_{1,n}^2 T^*(V_{0,1}) < T^*(u_0) - \tfrac{1}{2}\l_{1,n}^2T^*(V_{0,1}) < T^*(u_0)\ .$$
Notice that, since $u_n = NS(u_{0,n}) = NS(u(t_n))$, we may write $u(s_n) = u(t_n + (s_n - t_n)) = u_n(s_n - t_n)$ and hence, by (\ref{eq:ab}) we have (for a further subsequence in $n$)
\begin{equation}\label{eq:at}
u(s_n)  = \sum_{j=1}^{J} \frac{1}{\ljn} U_j \left( \frac{\cdot-\xjn}{\ljn}, \frac{s_n - t_n}{\ljn^2} \right) + w_n^{l,J}(s_n - t_n) + r_n^J(s_n - t_n)\ .
\end{equation}
Fix any $\e >0$.  Since $s_n - t_n = \frac{1}{2}T^*(V_{0,1})\l^2_{1,n} < T^*(V_{0,1})\l^2_{1,n}$ for all $n$, the orthogonality properties (\ref{eq:ad}) and (\ref{eq:ae}) and Young's inequality allow us to fix some $J$ so large that
\begin{equation}\label{eq:av}
\|r_n^J(s_n - t_n)\|_3 \lesssim \|r_n^J(s_n - t_n)\|_\hdhalf \leq \frac{\e}{4}\
\end{equation}
and
\begin{equation}\label{eq:aw}
\|w_n^{l,J}(s_n - t_n)\|_3 \lesssim \|w_n^{J}\|_3 \leq \frac{\e}{4}\
\end{equation}
for large enough $n$.  For such a $J$, we may take $n$ large enough that
\begin{equation}\label{eq:ax}
\max_{2\leq j \leq J} \left\|U_j\left(\frac{s_n - t_n}{\ljn^2}\right)\right\|_3  \leq \frac{\e}{2J}\ ,
\end{equation}
due to (\ref{eq:ar}) and the fact that $T^*(\voj) = +\infty$ for $j\geq 2$ which implies that $\lim_{t\to +\infty}\|U_j(t)\|_3 = 0$ by Theorem \ref{decay}.
We have therefore shown (due to (\ref{eq:at})) that for any $\e >0$, there exists $N\in \N$ such that
$$\left\|u(s_n) - \frac{1}{\l_{1,n}} U_1 \left( \frac{\cdot-x_{1,n}}{\l_{1,n}}, \frac{s_n - t_n}{\l_{1,n}^2} \right) \right\|_3 =
\left\|u(s_n) - \frac{1}{\l_{1,n}} U_1 \left( \frac{\cdot-x_{1,n}}{\l_{1,n}}, \tfrac{1}{2}T^*(V_{0,1}) \right) \right\|_3 < \e
$$
for all $n\geq N$.  A simple change of variables $y= \frac{x-x_{1,n}}{\l_{1,n}}$ shows that
$$\left\|\l_{1,n}u\left(\l_{1,n}\left[\ \cdot + \frac{x_{1,n}}{\l_{1,n}}\right],s_n\right) -  U_1 \left( \tfrac{1}{2}T^*(V_{0,1}) \right) \right\|_3 < \e$$
for $n\geq N$, and hence setting $\l_n = (\l_{1,n})^{-1}$ (note that $\l_n \to +\infty$ by (\ref{eq:aaaaa})) and $x_n = -{x_{1,n} / \l_{1,n}}$, we see that\footnote{As noted previously, for any $\eta \in (0,1)$ one can, of course, find a similar sequence which converges to ${U_1(\eta\cdot T^*(V_{0,1}))}$. }
$$\frac{1}{\l_n}u\left(\frac{\cdot - x_n}{\l_n}, s_n\right) \longrightarrow U_1(\tfrac{1}{2}T^*(V_{0,1}))$$
in $L^3$ as $n\to \infty$ and Theorem \ref{thm:b} is proved.  \hfill $\Box$
\subsection{Rigidity}
\noindent
In this section, we'll prove Theorem \ref{thm:c} which denies the possibility of the compactness of a critical element as established in Theorem \ref{thm:b}.  This, due to Theorems \ref{thm:a} and \ref{thm:b}, gives an alternative proof of the regularity criterion (\ref{eq:s}) in the case of mild solutions (which follows from the more general result in \cite{ess}), and concludes the critical element proof of Theorem \ref{thm:regularity}.
\\\\
We will use three key known results.  The first two concern systems of parabolic equations:  the backwards uniqueness theorem, Theorem \ref{thm:backwards}, proved in \cite{ess2,ess3} (see also \cite{ess}), and the following ``unique continuation" theorem, a proof of which can be found in \cite{ess} but which was already a well-known fact from the unique continuation theory of differential inequalities (see also \cite{ef}):
\begin{thm}[Unique Continuation]\label{thm:uniquecont}
Let $Q_{r,\d}:=  B_r(0)\times (-\delta,0)$ for some $r, \d > 0$, and suppose a vector-valued function $v$ and its distributional derivatives satisfy $v,\ v_t,\ \nabla v,\ \nabla^2 v \in L^2(Q_{r,\d})$ and there exist $c_0, C_k >0$ ($k\in \N$) such that
$|v_t - \Delta v| \leq c_0(|\nabla v| + |v|) \quad \textrm{a.e. on} \quad Q_{r,\d}$ and $|v(x,t)| \leq {C_k(|x|+\sqrt{-t})^k}$ for all $(x,t)\in Q_{r,\d}$.  Then $v(x,0)\equiv 0$ for all $x\in B_r(0)$.
\end{thm}
\noindent
The third result we will need is a local regularity criterion (with estimates) originating in \cite{caf} and generalized in \cite{necas} for the so-called ``suitable weak solutions" of the Navier-Stokes equations, which are defined as follows:
\begin{definition}[Suitable Weak Solutions]\label{def:suitable} Let $\O$ be an open subset of $\R^n$, $-\infty < T_1 < T_2 < +\infty$.  We call the pair $(u,p)$ a {\em suitable weak solution} of the Navier-Stokes equations in $\O \times (T_1,T_2)$ if
$u\in L^\infty((T_1,T_2);L^2(\O))\cap L^2((T_1,T_2);H^1(\O))$, $ p\in L^\infty((T_1,T_2);L^{{3/2}}(\O))$,
$u$ and $p$ satisfy NSE in distributions and the following local energy inequality holds:
$$\int_\O \varphi |u|^2\ dx + 2 \int_{T_1}^{t} \int_\O \varphi |\nabla u|^2\ dx\ ds \leq \qquad \qquad \qquad$$
$$\qquad \qquad \qquad \leq \int_{T_1}^{t} \int_\O \left\{  |u|^2(\Delta \varphi + \varphi_t) + u\cdot \nabla \varphi (|u|^2 + 2p) \right\}\ dx\ ds$$
for almost all $t\in (T_1,T_2)$ and for all non-negative cut-off functions $\varphi \in \mathcal{C}_0^\infty$ which vanish in a neighborhood of the parabolic boundary $(\O \times \{T_1\}) \cup (\partial \O \times [T_1,T_2])$.
\end{definition}
\noindent
We take the statement of the local regularity criterion from \cite{ess} (where a self-contained proof, based on the easier method in \cite{lin}, is presented).  It is the following:
\begin{lemma}[Local Smallness Regularity Criterion]\label{lemma:smallness}
There exist positive absolute constants $\e_0$ and $c_k$ for $k\in \N$ with the following property:  If a suitable weak solution $(u,p)$ of NSE on $Q_1$, where $Q_r:= B_r(0) \times (-r^2,0)$ for $r>0$, satisfies the condition
$$\int_{Q_1} (|u|^3 + |p|^{\frac{3}{2}})\ dx\ dt < \e_0\ ,$$
then $\nabla^{k-1} u$ is H\" older continuous on $\overline{Q_{\frac{1}{2}}}$ for any $k\in \N$, and for each $k$ we have the estimate
$$\max_{\overline{Q_\frac{1}{2}}}|\nabla^{k-1}u| \leq c_k\ .$$
\end{lemma}
\noindent
We will now prove Theorem \ref{thm:c} by using the backwards uniqueness, unique continuation and smallness criterion results (Theorem \ref{thm:backwards}, Theorem \ref{thm:uniquecont} and Lemma \ref{lemma:smallness}) to establish the following stronger version of Lemma \ref{lemma:zero}:
\begin{lemma}\label{lemma:final}
Suppose $u_0 \in \lt$ and set $u:=NS(u_0)$.  Suppose there exist some $T\in \R$ such that $0<T<+\infty$ and $T\leq T^*(u_0)$ and a sequence of numbers $\{s_n\}$ such that $0<s_n \nearrow T$ such that the following two properties hold:
$$
\begin{array}{l}
(1) \qquad \displaystyle{\sup_{t\in(0,T)}\|u(t)\|_3 < +\infty} \qquad \textrm{and} \\\\
(2) \qquad \displaystyle{\lim_{n\to \infty} \int_{|x|<R} |u(x,s_n)|^2\ dx = 0} \quad
\textrm{for any}\ R>0.
\end{array}
$$
Then $u_0 = u = 0$.
\end{lemma}
\noindent
(The case $T<T^*(u_0)$ implies Lemma \ref{lemma:zero} due to the continuity properties of mild solutions.) Note the immediate corollary that there exists no $u_0 \in \lt$ such that $T^*(u_0) < +\infty$ and (1) and (2) hold with $T=T^*(u_0)$, since the conclusion of Lemma \ref{lemma:final} implies that $T^*(u_0) = T^*(0) = +\infty$ by Theorem \ref{unique}.  This will exactly rule out the critical element produced in Theorem \ref{thm:a} and prove the regularity criterion (\ref{eq:s}).  Let us postpone the proof of Lemma \ref{lemma:final} for the moment.
\\\\
\textbf{Proof of Theorem \ref{thm:c}:}
\quad
As in Theorem \ref{thm:b}, define $A_c$ by (\ref{eq:an}) and suppose $A_c<+\infty$, and let $u_{0} \in \hdhalf$ with associated mild solution $u:=NS(u_0)$ satisfy $T^*:= T^*(u_{0}) < +\infty$ and $\|u\|_{L^\infty((0,T^*(u_{0}));\hdhalf)}=A_c$.
\\\\
Fix any $\{t_n\} \subset [0,T^*)$ such that $t_n \nearrow T^*$, and let $s_n\nearrow T^*$, $\l_n \to +\infty$ and $\{x_n\} \subset \rt$ be the associated sequences guaranteed by Theorem \ref{thm:b}. Using the conclusion of the same theorem, define
$$v_n(x) = \frac{1}{\l_n}u\left(\frac{x-x_n}{\l_n},s_n\right)$$
for any $x\in \rt$ and pass to a subsequence in $n$ so that $v_n \to \bar{v}$ in $\lt$ for some $\bar{v}\in \lt$.  We now make the following important claim:
\begin{clm}\label{clm:zero} For $s_n \nearrow T^*$ as above and for any $R>0$, $$\lim_{n\to \infty} \int_{B_R(0)} |u(x,s_n)|^2\ dx = 0\ .$$
\end{clm}
\noindent
\textbf{Proof of Claim \ref{clm:zero}:}
\quad Fix $R>0$.  We calculate:
$$\int_{|x| \leq R} |u(x,s_n)|^2\ dx = \int_{|x| \leq R} |\l_nv_n(\l_n x+x_n)|^2\ dx =\qquad \qquad  $$$$\qquad \qquad =(\l_n)^{-1} \int_{|y-x_n|\leq \l_nR}|v_n(y)|^2\ dy
=:(\l_n)^{-1} \|v_n\|^2_{L^2(B_n)}\ .$$
For any small $\e >0$, define $B_{\e \l_n R} := B_{\e \l_n R} (0)$.  We then estimate
$$(\l_n)^{-1} \|v_n\|^2_{L^2(B_n)} = (\l_n)^{-1} \|v_n\|^2_{L^2(B_n\cap B_{\e \l_n R})} + (\l_n)^{-1} \|v_n\|^2_{L^2(B_n\cap B_{\e \l_n R}^c)}\qquad $$
$$\quad \qquad \qquad \qquad \qquad \leq(\l_n)^{-1} \|v_n\|^2_{L^3(\rt)}|B_{\e \l_n R}|^\frac{1}{3} + (\l_n)^{-1} \|v_n\|^2_{L^3(B_n\cap B_{\e \l_n R}^c)}|B_n|^\frac{1}{3} $$
$$\leq C(A_c)^2\e R + 4CR \|\bar{v}\|^2_{L^3(B_{\e \l_n R}^c)} \qquad $$
for sufficiently large $n$.  The first term of the last line is arbitrarily small for small $\e$, and for fixed $\e$ the second term is arbitrarily small for sufficiently large $n$ due to the fact that $\l_n \to +\infty$.  This proves the claim. \hfill $\Box$
\\\\
Assuming Lemma \ref{lemma:final} holds and taking $T=T^*$, Claim \ref{clm:zero} clearly completes the proof of Theorem \ref{thm:c}.  Therefore it remains only to establish the following:
\\\\
\textbf{Proof of Lemma \ref{lemma:final}:} \quad Since $T=T^*=T^*(u_0)$ is the application we seek, we will consider only that case (the case $T<T^*(u_0)$ is even easier).  Assume therefore that $u_0\in \lt$ is such that $T^*(u_0) < +\infty$, $$\sup_{t\in (0,T^*(u_0))} \|u(t)\|_3 < +\infty$$ where $u=NS(u_0)$, and there exists $\{s_n\}$ with $s_n\nearrow T^*$ as $n\to \infty$ such that $$\lim_{n\to \infty} \int_{|x|<R} |u(x,s_n)|^2\ dx = 0$$ for any fixed $R>0$.
\\\\
Recall that $\P F = F - \nabla \frac{1}{\D} \nabla \cdot F$ for a vector field $F$, and for a tensor $G$ we can write
$-  \frac{1}{\D} (\nabla \cdot (\nabla \cdot G)) =   (\frac{\nabla }{\sqrt{-\D}}\cdot (\frac{\nabla }{\sqrt{-\D}} \cdot G)) = \RR_i \RR_j G_{ij}$, where $(\RR_j g)^{\wedge} (\xi) = \frac{i \xi_j}{|\xi|}\hat g (\xi)$ and we sum over repeated indices.  Therefore, letting $G = u\otimes u$, $F=\nabla \cdot G$ and defining $p:=\RR_i \RR_j u_i u_j$, we have
$$ \P \nabla \cdot (u\otimes u)  =  \nabla \cdot (u\otimes u) + \nabla p  \ .$$
Note that $p$ is well-defined as a distribution since $\RR_k : L^{\frac{3}{2}} \to L^{\frac{3}{2}}$ boundedly for each $k$ by the Calderon-Zygmund theory, so that
\begin{equation}\label{eq:aae}
\|p(t)\|_\frac{3}{2} \leq C\|u(t)\|_3^2\ ,
\end{equation}
hence\footnote{We adopt the notation of \cite{ess}, setting $L_{p,\infty} = L^\infty_t L^p_x$ for $p\geq 1$, and denote the corresponding norm by $\|\cdot \|_{p,\infty}$.} $p\in L_{\frac{3}{2},\infty}$.  Moreover, $p$ satisfies the equation
$-\D p(t) = \partial_i \partial_j u_i(t) u_j(t)$
in distributions and therefore $p(t)$ is smooth in $x$ for any $t\in (0,T^*)$ since $u(t)$ is smooth for such t.  Moreover, setting $B_r = B_r(x_0)$ for any fixed $x_0 \in \rt$ and any $r>0$, we have the following estimates from the classical elliptic theory (see, e.g., \cite{necas})
\begin{equation}\label{eq:aad}
\|p(\cdot, t)\|_{C^{k,\alpha}(B_r)} \leq c_k (\sum_{i,j} \|u_i(\cdot, t) u_j(\cdot, t)\|_{C^{k,\alpha}(B_{2r})} + \|p(\cdot, t)\|_{L^{\frac{3}{2}}(\rt)} )\ .
\end{equation}
For any $t\in (0, T^*)$, we now have
$$
u(t) = e^{t\Delta }u_0 - \int_{0}^t e^{(t-s)\Delta} [\nabla \cdot (u\otimes u) + \nabla p](s) \ ds\  .
$$
It is clear that $u\in \cap_{0<T<T^*}L^2_{\textrm{loc unif},x}L^2_t(\rt \times(0,T))$, due to the fact that $u\in L_{3,\infty}$.  Therefore $u$ satisfies (see \cite{lr}, Theorem 11.2)
\begin{equation}\label{eq:aaa}
u_t - \Delta u + \nabla \cdot (u\otimes u) + \nabla p = 0\ , \qquad \nabla \cdot u = 0
\end{equation}
in distributions.  We now claim that the pair $(u,p)$ locally forms a suitable weak solution, in the sense of Definition \ref{def:suitable}.
\\\\
Recall that, by construction, $u\in  L^5( \rt \times (0,T))$ for any $T<T^*$, so that the ``Ladyzhenskaya-Prodi-Serrin" regularity condition (see, e.g., \cite{ess}) implies that $u$ is smooth on $\rt \times (0,T^*)$ (or, alternatively, the construction gives $u\in L^\infty(\rt \times (\e,T^*-\e)$ for any $\e>0$ and hence is smooth, see, e.g., \cite{lr}, Prop. 15.1).  Therefore (\ref{eq:aae}) and (\ref{eq:aad}) together imply that $p(t)$ is smooth in $x$ for $0<t<T^*$ and $\nabla^k_x p \in L^\infty_\textrm{loc}(\rt \times (0,T^*))$ for any $k$.
\\\\
We clearly have $u\in L^\infty((0,T^*);L^2(\O))$ for any bounded $\O \subset \subset \rt$.  Since $u$ and $p$ are sufficiently smooth
on $\rt \times (0,T^*)$ we may multiply the equation (\ref{eq:aaa}) by $u\varphi$ for $\varphi$ as in Definition \ref{def:suitable} and integrate by parts to derive the local energy inequality in $\tilde{\O} \times (0,T^*)$ for a suitably larger region such that $\O \subset \subset \tilde{\O}$, which implies that $u\in L^2 ((0,T^*);H^1(\O))$.  (Indeed, the integrations by parts happen only at the level of integration in $x$, and then $(u,p)\in L_{3,\infty} \times L_{\frac{3}{2},\infty}$ allows integration in $t$.)  Noting also that we have
\begin{equation}\label{eq:aab}
\int_0^{T^*} \int_{\rt} \left( |u|^3 + |p|^\frac{3}{2} \right) \ dx\ dt \leq T^* \left(\|u\|_{3,\infty}^3 + \|p\|_{\frac{3}{2},\infty}^{\frac{3}{2}}\right) < \infty\ ,
\end{equation}
we see that $p\in L^{\frac{3}{2}}(\O \times (0,T^*))$.  The pair $(u,p)$ now clearly satisfies all the requirements to form a suitable weak solution to NSE in $\O \times (0,T^*)$ for any $\O \subset \subset \rt$.
\\\\
Equation (\ref{eq:aab}) implies moreover that for any $\e >0$ one can find $R>0$ large enough that
$$\int_0^{T^*}\int_{B_{\sqrt{T^*}}(x_0)} \left(|u|^3 + |p|^{\frac{3}{2}}\right)\ dx\ dt \leq \e \ $$
for all $x_0$ such that $|x_0|\geq R$.  Therefore, due to Lemma \ref{lemma:smallness} (and an appropriate scaling argument), we see that $u$ is smooth on $\Omega:=(\rt \backslash B_{R_0}(0))\times [({3/4})T^*, T^*]$ for some sufficiently large $R_0$, and by shifting spatial regions, we have the global bounds
\begin{equation}\label{eq:aac}
\max_{\overline{\Omega}}|\nabla^{k-1}u| \leq c_k\ , \quad k=1,2,3,\dots
\end{equation}
Since $u$ is continuous up to $T^*$ outside $B_{R_0}$, Claim \ref{clm:zero} now implies that $u(x,T^*) \equiv 0$ for all $x$ such that  $|x|\geq R_0$.
We therefore also have $\omega:= \nabla_x \times u \equiv 0$ on $(\rt \backslash B_{R_0}) \times \{T^*\}$. Taking the curl of equation (\ref{eq:aaa}), we see that the vorticity $\omega$ satisfies the inequality
$$|\omega_t -\D \omega| = |(u\cdot \nabla) \omega - (\omega \cdot \nabla) u| \leq c(|\omega| + |\nabla \omega|)$$
in $\Omega$ for some $c>0$ due to (\ref{eq:aac}).  Since $\omega$ is bounded and smooth in $\Omega$, we can apply the backwards uniqueness theorem, Theorem \ref{thm:backwards}, to conclude that $\omega \equiv 0$ in $\Omega$.
\\\\
Now, in $\tilde{\Omega}:=\rt \times (({3/4})T^*, ({7/8})T^*]$ we have $u$, and hence $\omega$, is smooth, hence $D^m_t D^k_x \omega \in L^2_{\textrm{loc}}(\tilde{\Omega})$ for all $m,k \geq 0$.  We can therefore apply the theorem of unique continuation, Theorem \ref{thm:uniquecont}, to $\omega$, since $\omega \equiv 0$ in $\tilde{\Omega}\cap \Omega$, to conclude (by appropriate shifting of local regions) that $\omega \equiv 0$ in $\tilde{\Omega}$.  Therefore, due to the divergence-free condition, $u$ is harmonic in $\tilde{\Omega}$ which, along with the fact that $u\in L_{3,\infty}$, implies $u\equiv 0$ in $\tilde{\Omega}$.  It is then easy to apply Theorem \ref{thm:backwards} again (see the proof of Lemma \ref{lemma:zero}) to see that $u_0=u = 0$  which proves Lemma \ref{lemma:final} and concludes the proof of Theorem \ref{thm:c} and Theorem \ref{thm:regularity}. \hfill $\Box$


\end{document}